\documentclass[12pt]{article}
\usepackage{amsmath,amsfonts,amssymb,amsthm}
\input amssym.def

\begin{document}
\def \Z{\mathbb Z}
\def \C{\mathbb C}
\def \R{\mathbb R}
\def \Q{\mathbb Q}
\def \N{\mathbb N}

\def \A{{\mathcal{A}}}
\def \D{{\mathcal{D}}}
\def \E{{\mathcal{E}}}
\def \E{{\mathcal{E}}}
\def \H{\mathcal{H}}
\def \S{{\mathcal{S}}}
\def \U{{\mathcal{U}}}
\def \Y{{\mathcal{Y}}}

\def \wt{{\rm wt}}
\def \tr{{\rm tr}}
\def \add{{\rm add}}
\def \span{{\rm span}}
\def \Res{{\rm Res}}
\def \Der{{\rm Der}}
\def \End{{\rm End}}
\def \Ind {{\rm Ind}}
\def \Irr {{\rm Irr}}
\def \Aut{{\rm Aut}}
\def \GL{{\rm GL}}
\def \Hom{{\rm Hom}}
\def \mod{{\rm mod}}
\def \ann{{\rm Ann}}
\def \ad{{\rm ad}}
\def \rank{{\rm rank}\;}
\def \<{\langle}
\def \>{\rangle}

\def \g{{\mathfrak{g}}}
\def \h{{\hbar}}
\def \k{{\mathfrak{k}}}
\def \sl{{\mathfrak{sl}}}
\def \gl{{\mathfrak{gl}}}

\def \be{\begin{equation}\label}
\def \ee{\end{equation}}
\def \bex{\begin{example}\label}
\def \eex{\end{example}}
\def \bl{\begin{lem}\label}
\def \el{\end{lem}}
\def \bt{\begin{thm}\label}
\def \et{\end{thm}}
\def \bp{\begin{prop}\label}
\def \ep{\end{prop}}
\def \br{\begin{rem}\label}
\def \er{\end{rem}}
\def \bc{\begin{coro}\label}
\def \ec{\end{coro}}
\def \bd{\begin{de}\label}
\def \ed{\end{de}}

\newcommand{\q}{\bf q}
\newcommand{\m}{\bf m}
\newcommand{\n}{\bf n}
\newcommand{\nno}{\nonumber}
\newcommand{\nord}{\mbox{\scriptsize ${\circ\atop\circ}$}}
\newtheorem{thm}{Theorem}[section]
\newtheorem{prop}[thm]{Proposition}
\newtheorem{coro}[thm]{Corollary}
\newtheorem{conj}[thm]{Conjecture}
\newtheorem{example}[thm]{Example}
\newtheorem{lem}[thm]{Lemma}
\newtheorem{rem}[thm]{Remark}
\newtheorem{de}[thm]{Definition}
\newtheorem{hy}[thm]{Hypothesis}
\makeatletter \@addtoreset{equation}{section}
\def\theequation{\thesection.\arabic{equation}}
\makeatother \makeatletter

\begin{center}{\Large \bf Twisted regular representations of
vertex operator algebras}
\end{center}

\begin{center}
{Haisheng Li \\
Department of Mathematical Sciences\\
Rutgers University, Camden, NJ 08102\\
Jiancai Sun\footnote{Partially supported by China NSF grants}\\
Department of Mathematics\\
Shanghai University, Shanghai 200444, China}
\end{center}

\begin{abstract}
This paper is to study what we call twisted regular representations for vertex operator algebras. 
Let $V$ be a vertex operator algebra, let $\sigma_1,\sigma_2$ be commuting finite-order automorphisms of $V$
and let $\sigma=(\sigma_1\sigma_2)^{-1}$.
Among the main results, for any $\sigma$-twisted $V$-module $W$ and any nonzero complex number $z$, 
we construct a weak $\sigma_1\otimes \sigma_2$-twisted $V\otimes V$-module
$\mathfrak{D}_{\sigma_1,\sigma_2}^{(z)}(W)$ inside $W^{*}$. 
Let $W_1,W_2$ be $\sigma_1$-twisted, $\sigma_2$-twisted $V$-modules, respectively. 
We show that $P(z)$-intertwining maps from $W_1\otimes W_2$ to $W^{*}$
are the same as homomorphisms of weak $\sigma_1\otimes \sigma_2$-twisted $V\otimes V$-modules 
from $W_1\otimes W_2$ into $\mathfrak{D}_{\sigma_1,\sigma_2}^{(z)}(W)$.
We also show that a $P(z)$-intertwining map from $W_1\otimes W_2$ to $W^{*}$ is equivalent to
an intertwining operator of type $\binom{W'}{W_1\; W_2}$, which is a twisted version of a result of Huang and Lepowsky.
Finally, we show that for each $\tau$-twisted $V$-module $M$ with $\tau$ any finite-order automorphism of $V$, 
the coefficients of the $q$-graded trace function lie in 
$\mathfrak{D}_{\tau,\tau^{-1}}^{(-1)}(V)$, which generate a $\tau\otimes \tau^{-1}$-twisted $V\otimes V$-submodule isomorphic to
$M\otimes M'$.
\end{abstract}

\section{Introduction}
In \cite{li-reg}, a theory of what were called regular representations of 
vertex operator algebras was developed and among the main results 
an analogue of the classical Peter-Weyl theorem (cf. \cite{bd}) was established. More specifically,
 let $V$ be a vertex operator algebra.  
For any $V$-module $W$ and any nonzero complex number $z$, we constructed a (weak) $V\otimes V$-module 
${\mathfrak{D}}_{P(z)}(W)$ inside the full dual space $W^{*}$. It was proved that for any $V$-modules 
$W_{1}$ and $W_{2}$, an intertwining operator of type ${W'\choose W_{1}W_{2}}$ 
amounts to a $V\otimes V$-module homomorphism from $W_{1}\otimes W_{2}$ into ${\mathfrak{D}}_{P(z)}(W)$.
Then it follows that the socle of ${\mathfrak{D}}_{P(z)}(W)$, i.e., the sum of irreducible $V\otimes V$-submodules, 
is decomposed canonically into irreducible submodules with the spaces of intertwining operators as the multiplicity spaces. 
In the special case with $W=V$, we obtained a Peter-Weyl type decomposition of the socle of ${\mathfrak{D}}_{P(z)}(V)$.

 In the theory of vertex operator algebras, intertwining operator, fusion rule, and contragredient dual, 
 which were introduced by Frenkel-Huang-Lepowsky (see \cite{fhl}), are of fundamental importance. 
In developing their tensor product theory, Huang and Lepowsky introduced a notion of $P(z)$-intertwining map  (see \cite{hl2})
where $z$ is a nonzero complex number. For a triple of $V$-modules $W_1,W_2$, and $W$, 
a $P(z)$-intertwining map of type $\binom{W'}{W_1\;W_2}$ by definition is a linear map from $W_1\otimes W_2$ to $W^{*}$, 
which satisfies a Jacobi-type identity.
The fact is that any intertwining operator $\mathcal{Y}(\cdot,x)$ of type
$\binom{W'}{W_1\;W_2}$ with the formal variable $x$ evaluated at $z$ becomes a $P(z)$-intertwining map.
It was proved therein that the evaluation $x=z$ gives rise to a linear isomorphism from the space of intertwining operators 
to the space of $P(z)$-intertwining maps. Just as in Huang-Lepowsky's tensor product theory, 
$P(z)$-intertwining map plays an important role in the study of regular representations. 
It was shown in \cite{li-reg} that a $P(z)$-intertwining map of type ${W'\choose W_{1}W_{2}}$ 
is the same as a $V\otimes V$-module homomorphism from $W_{1}\otimes W_{2}$ into ${\mathfrak{D}}_{P(z)}(W)$.

Needless to say, regular representation is closely related to Huang-Lepowsky's tensor product theory. 
(Relations between them were further studied in \cite{Li-reg-hl}.)
In addition, regular representation also has intrinsic connections with Zhu's $A(V)$-theory and its generalizations,
 which enable us to obtain new (arguably shorter) proofs for several known theorems. 
Note that the Zhu algebra $A(V)$ of a vertex operator algebra $V$ is a quotient space of $V$ (see \cite{zhu, zhu1}), while
the Frenkel-Zhu $A(V)$-bimodule $A(W)$ associated to a $V$-module $W$ is a quotient space of $W$ (see \cite{fz}).
It was proved in \cite{li-reg-Zhu} that for any $V$-module $W$, $A(W)^{*}$ viewed as a subspace of $W^{*}$ 
 coincides with the vacuum space of ${\mathfrak{D}}_{P(-1)}(W)$, which is naturally an $A(V)\otimes A(V)$-module
 and an $A_n(V)$-bimodule through a canonical anti-automorphism of $A_n(V)$. 
 For any irreducible $A(V)$-module $U$, by making use of the weak $V\otimes V$-module ${\mathfrak{D}}_{P(-1)}(V)$, 
 we can straightforwardly obtain an irreducible $V$-module with $U$ as its vacuum space.
On the other hand, using the result that an intertwining operator of type ${W'\choose W_{1}W_{2}}$ 
amounts to a $V\otimes V$-module homomorphism from $W_{1}\otimes W_{2}$ into ${\mathfrak{D}}_{P(z)}(W)$,
we obtained (loc cit) a new (shorter) proof of Frenkel-Zhu's fusion rule theorem (see \cite{fz}; cf. \cite{li-thesis, li-fusion}). 
The relation between the $A_{n}(V)$-theory (see \cite{dlm-anv}) and regular representation
was also studied (see \cite{Li-reg-anv}).
In a recent work \cite{Li-m-n-bimodule}, regular representation was used to reexamine the $A_n(V)$-$A_m(V)$ bimodules
 introduced by Dong and Jiang (see \cite{dj-1}) and the main results therein were recovered.

In the representation theory of vertex operator algebras, twisted module, which originated in the construction of the celebrated 
moonshine module  (see \cite{flm}),  plays a significant role.
Contragedient duals of twisted modules and intertwining operators for triples of twisted modules were studied by Xu (see \cite{xu}).
More specifically, let $\sigma_1,\sigma_2,\sigma_3$ be mutually commuting finite-order automorphisms of a vertex operator algebra $V$
and let $W_r$ be a $\sigma_r$-twisted $V$-module for $r=1,2,3$. Then Xu defined a notion of intertwining operator of type
${W_3\choose W_1\;W_2}$, and he showed that if there exists a nonzero intertwining operator 
of this type, then a necessary condition is that $\sigma_3=\sigma_1\sigma_2$. 
In view of the various applications of regular representations, it is of conceptual and practical importance to study
 regular representations associated to twisted modules.

In this paper, we study twisted regular representations of vertex operator algebras and we
present the basic results to lay the foundation for subsequent studies.
Let $V$ be a vertex operator algebra, let $\sigma_1,\sigma_2$ be commuting finite-order automorphisms of $V$, 
and set $\sigma=(\sigma_1\sigma_2)^{-1}$.
Among the main results, for any given $\sigma$-twisted $V$-module $W$ and nonzero complex number $z$, 
we construct a (weak) $\sigma_1\otimes \sigma_2$-twisted $V\otimes V$-module
$\mathfrak{D}_{\sigma_1,\sigma_2}^{(z)}(W)$ inside $W^{*}$. 
Let $W_1,W_2$ be $\sigma_1$-twisted, $\sigma_2$-twisted $V$-modules, respectively. 
We show that $P(z)$-intertwining maps from $W_1\otimes W_2$ to $W^{*}$
are the same as homomorphisms of weak $\sigma_1\otimes \sigma_2$-twisted $V\otimes V$-modules 
from $W_1\otimes W_2$ into $\mathfrak{D}_{\sigma_1,\sigma_2}^{(z)}(W)$.
We also show that a $P(z)$-intertwining map from $W_1\otimes W_2$ to $W^{*}$ exactly amounts to
an intertwining operator of type $\binom{W'}{W_1\; W_2}$, which slightly generalizes a result of Huang and Lepowsky 
for (untwisted) modules. We also show that for each $\tau$-twisted $V$-module $M$ with $\tau$ any finite-order automorphism of $V$, 
the coefficients of the $q$-graded trace function lie in 
$\mathfrak{D}_{\tau,\tau^{-1}}^{(-1)}(V)$, which generate a $\tau\otimes \tau^{-1}$-twisted $V\otimes V$-submodule isomorphic to
$M\otimes M'$. This gives a direct (canonical) association of irreducible twisted modules to ``characters.''

This paper is organized as follows: In Section 2, we recall basic notions and results such as twisted modules, contragredient dual.
Section 3 is the core, in which we present the construction of weak $\sigma_1\otimes \sigma_2$-twisted $V\otimes V$-module 
$\mathfrak{D}_{\sigma_1,\sigma_2}^{(z)}(W)$. In Section 4, we study the interplay of
intertwining operators, $P(z)$-intertwining maps, and $\mathfrak{D}_{\sigma_1,\sigma_2}^{(z)}(W)$. 
We also present a connection of graded trace functions for $\tau$-twisted $V$-modules with
$\mathfrak{D}_{\tau,\tau^{-1}}^{(-1)}(V)$.

For this paper, we use $\N$ for the set of nonnegative integers and $\Z_{+}$ for the set of positive integers,
in addition to the standard notations $\Z, \Q, \R$, and $\C$ for the sets of integers, rational numbers, real numbers, and complex numbers,
respectively.

\section{Preliminaries}
In this section, we recall mostly from \cite{flm} and \cite{fhl} some basic notions and results, including contragredient module, 
twisted module, and right twisted module.

We begin by recalling the definition of a vertex algebra.
 A {\em  vertex algebra } is a vector space $V$ equipped with a linear map
\begin{eqnarray*}
Y(\cdot,x):&& V\rightarrow (\End V)[[x,x^{-1}]]\\
&&v\mapsto Y(v,x),
\end{eqnarray*}
satisfying the condition that 
\begin{eqnarray}
Y(u,x)v\in V((x))\   \   \   \   \mbox{ for }u,v \in V,
\end{eqnarray}
$Y({\bf 1},x)=1$, 
$$Y(u,x){\bf 1}\in V[[x]] \  \ \text{ and }\ \ \lim_{x\rightarrow 0}Y(u,x){\bf 1}=u\quad \text{for }u\in V,$$
and 
\begin{eqnarray}
& &x_{0}^{-1}\delta\left(\frac{x_{1}-x_{2}}{x_{0}}\right)Y(u,x_{1})Y(v,x_{2})
-x_{0}^{-1}\delta\left(\frac{x_{2}-x_{1}}{-x_{0}}\right)Y(v,x_{2})Y(u,x_{1})   \nonumber\\
&&\hspace{2cm}   =x_{2}^{-1}\delta\left(\frac{x_{1}-x_{0}}{x_{2}}\right)Y(Y(u,x_{0})v,x_{2})
\end{eqnarray}
(the {\em Jacobi identity}) for $u,v\in V$.

A {\em conformal vertex algebra} is a vertex algebra $V$ equipped with a distinguished vector $\omega$, 
called the {\em conformal vector}, such that  
\begin{eqnarray}\label{relation-vir}
[L(m),L(n)]=(m-n)L(m+n)+\frac{1}{12}(m^3-m)\delta_{m+n,0}c_V
\end{eqnarray}
for $m,n\in \Z$, where $Y(\omega,x)=\sum_{n\in \Z}L(n)x^{-n-2}$ and $c_V\in \C$, called the {\em central charge,}
and such that 
\begin{eqnarray}\label{L(-1)-cva}
Y(L(-1)v,x)=\frac{d}{dx}Y(v,x)\quad \text{ for }v\in V,
\end{eqnarray}
\begin{eqnarray}
 V=\oplus_{n\in \Z}V_{(n)}, \quad \text{where }V_{(n)}=\{ v\in V\  |\  L(0)v=nv\}.
\end{eqnarray}
In addition, it is assumed that $L(1)$ is locally nilpotent on $V$. 
 
If $V$ is a conformal vertex algebra, we have
\begin{eqnarray}
{\bf 1}\in V_{(0)},\quad \omega\in V_{(2)}.
\end{eqnarray} 

 Note that for any vertex algebra $V$,  we have a canonical derivation $\D$ which is defined by $\D (v)=v_{-2}{\bf 1}$ for $v\in V$. 
 Furthermore, the following relations hold:
 \begin{eqnarray}\label{D-va}
 [\D, Y(v,x)]=Y(\D (v),x)=\frac{d}{dx}Y(v,x)
 \end{eqnarray}
 for $v\in V$.  If $V$ is a conformal vertex algebra, it follows that $L(-1)=\D$.
 
A {\em vertex operator algebra} is a conformal vertex algebra $V$ satisfying the condition that
\begin{eqnarray*}
&&\dim V_{(n)}<\infty\quad \text{ for all }n\in \Z,\\
&&V_{(n)}=0\quad \text{  for $n$ sufficiently negative.}
\end{eqnarray*} 
 
Recall that an {\em automorphism} of a vertex algebra $V$ is an invertible linear operator $\sigma$ on $V$ such that
\begin{eqnarray*}
\sigma({\bf 1})={\bf 1},\    \    \   \sigma Y(u,x)v=Y(\sigma (u),x)\sigma (v)\quad  \mbox{ for }u,v\in V. 
\end{eqnarray*}
If $V$ is a conformal vertex algebra, an automorphism of $V$ by definition preserves the conformal vector of $V$, so that  
it preserves each weight subspace of $V$.

Let $V$ be a vertex algebra, $\sigma$ a finite-order automorphism, and 
let $T$ be a positive integer such that $\sigma^T=1$.
For $r\in \Z$, set
\begin{eqnarray}
V^{r}=\left\{v\in V\ |\  \sigma(v)=e^{{2r\pi i\over T}}v\right\}.
\end{eqnarray}
Note that for $r,s\in \Z$, if $r\equiv s \  (\mod \ T)$, then $V^{r}=V^{s}$.
We have
\begin{eqnarray}
V=V^{0}\oplus V^{1}\oplus\cdots \oplus V^{T-1}.
\end{eqnarray}

With this setting,  we have (see  \cite{flm}, \cite{ffr}, \cite{dong}):

\bd{dtwistedmodule}
{\em A {\em $\sigma$-twisted $V$-module}  is a vector space $W$ equipped with a linear map 
\begin{eqnarray*}
Y_{W}(\cdot,x):&&V\rightarrow ({\rm End}W)[[x^{\frac{1}{T}},x^{-\frac{1}{T}}]],\\
&&v\mapsto Y_{W}(v,x),
\end{eqnarray*}
 satisfying the following conditions that
$Y_{W}({\bf 1},x)=1,$
\begin{eqnarray}
Y_{W}(v,x)w\in x^{-\frac{j}{T}}W((x))\   \    \mbox{ for }
v\in V^{j},\  w\in W, \ j\in \Z,
\end{eqnarray}
and that for $u\in V^{j}, \  v\in V$ with $j\in \Z$,
\begin{eqnarray}\label{twistedj}
& &x_{0}^{-1}\delta\left(\frac{x_{1}-x_{2}}{x_{0}}\right)Y_{W}(u,x_{1})Y_{W}(v,x_{2})
-x_{0}^{-1}\delta\left(\frac{x_{2}-x_{1}}{-x_{0}}\right)Y_{W}(v,x_{2})Y_{W}(u,x_{1})
\nonumber\\
& &\hspace{1.5cm}=
x_{2}^{-1}\delta\left(\frac{x_{1}-x_{0}}{x_{2}}\right)
\left(\frac{x_{1}-x_{0}}{x_{2}}\right)^{-\frac{j}{T}}
Y_{W}(Y(u,x_{0})v,x_{2})
\end{eqnarray}
(the {\em $\sigma$-twisted Jacobi identity}).}
\ed

Note that the notion of $\sigma$-twisted $V$-module for $\sigma=1$ simply reduces to
the notion of $V$-module.

For any $\sigma$-twisted $V$-module $(W,Y_W)$, we have (cf. \cite{dlm1}):
\begin{eqnarray}
Y_{W}(\mathcal{D}v,x)={d\over  dx}Y_{W}(v,x)\    \    \   \mbox{  for }v\in V.
\end{eqnarray}

The following is a straightforward fact (cf. \cite{dlm-twisted}):

\bl{lem-sigma-tau}
Suppose that $\sigma$ and $\tau$ are commuting automorphisms of $V$.
Then for any $\sigma$-twisted $V$-module $(W,Y_W)$, $(W,Y_W\circ \tau)$
is also a $\sigma$-twisted $V$-module, where $(Y_W\circ \tau)(v,x)=Y_{W}(\tau v,x)$ for $v\in V$.
\el

Let $V$ be a conformal vertex algebra. A $\sigma$-twisted module for $V$ viewed as a vertex algebra 
is called a {\em weak $\sigma$-twisted $V$-module.}
If $V$ is a vertex operator algebra, a {\em $\sigma$-twisted $V$-module} is a weak $\sigma$-twisted $V$-module 
$W$ satisfying the condition that 
$$W=\oplus_{h\in \C}W_{(h)}, \  \  \text{where }W_{(h)}=\{ w\in W\ |\ L(0)w=hw\},$$
such that $\dim W_{(h)}<\infty$ for all $h\in \C$ and for every fixed number $h\in \C$, 
$W_{(h+n)}=0$ for $n\in \frac{1}{T}\Z$ sufficiently negative.

Assume that $V$ is a vertex operator algebra and $(W,Y_{W})$ is a $\sigma$-twisted $V$-module.  
Just as in \cite{fhl} for a $V$-module, set
$$W'=\oplus_{h\in \C}W_{(h)}^{*},$$
and for $v\in V$,  define
$$Y'_{W}(v,x)\in (\End W')[[x^{\frac{1}{T}},x^{-\frac{1}{T}}]]$$
 by
\begin{eqnarray}
\<Y'_{W}(v,x)w',w\>=\<w',Y_{W}(e^{xL(1)}(-x^{-2})^{L(0)}v,x^{-1})w\>
\end{eqnarray}
for $w'\in W',\  w\in W$. 
The following result, due to Xu (see \cite{xu}), is a twisted analogue of a well-known theorem of 
Frenkel-Huang-Lepowsky (see \cite{fhl}): 

\bp{ptdual}
Let $(W,Y_{W})$ be a $\sigma$-twisted $V$-module. Then $(W',Y_{W}')$ carries 
the structure of a $\sigma^{-1}$-twisted $V$-module. Furthermore, $W$ is irreducible if and only if $W'$ is irreducible.
\ep

\br{rem-sigma-inverse}
{\em In view of Proposition \ref{ptdual}, if $\{ W_{\alpha}\ |\ \alpha \in S\}$ is a complete set of equivalence class representatives of 
irreducible $\sigma$-twisted $V$-modules, then
 $\{ W_{\alpha}'\ |\ \alpha \in S\}$ is a complete set of equivalence class representatives 
of irreducible $\sigma^{-1}$-twisted $V$-modules.
 It also follows from Proposition \ref{ptdual} that the category of $\sigma$-twisted $V$-modules is semisimple 
if and only if the category of $\sigma^{-1}$-twisted $V$-modules is semisimple. }
 \er

Let $V$ be a vertex operator algebra, $(W,Y_W)$ a $\sigma$-twisted $V$-module. For $v\in V$, set
\begin{eqnarray}
Y_{W}^{o}(v,x)=Y_{W}\left(e^{xL(1)}(-x^{-2})^{L(0)}v,x^{-1}\right)\in(\End
W)[[x^{\frac{1}{T}},x^{-\frac{1}{T}}]].
\end{eqnarray}
Note that 
\begin{eqnarray}\label{Yo-twisted}
Y_{W}^{o}(v,x)w\in W((x^{-\frac{1}{T}}))\quad \text{ for }v\in V, \  w\in W
\end{eqnarray}
and that the following opposite twisted Jacobi identity holds on $W$ for $u\in V^j,\ v\in V$: 
\begin{eqnarray}\label{opposite-Jacobi}
& &x_{0}^{-1}\delta\left(\frac{x_{1}-x_{2}}{x_{0}}\right)
Y_{W}^{o}(v,x_{2})Y_{W}^{o}(u,x_{1}) 
-x_{0}^{-1}\delta\left(\frac{x_{2}-x_{1}}{-x_{0}}\right)
Y_{W}^{o}(u,x_{1})Y_{W}^{o}(v,x_{2})\  \  \  \   \nonumber\\
&&\hspace{2cm}   =x_{2}^{-1}\delta\left(\frac{x_{1}-x_{0}}{x_{2}}\right)\left(\frac{x_{1}-x_{0}}{x_{2}}\right)^{\frac{j}{T}}
Y_{W}^{o}(Y(u,x_{0})v,x_{2}).
\end{eqnarray}

With this, we formulate a twisted analogue of the notion of right module (cf. \cite{hl1}).

\bd{dright-module}
{\em Let $V$ be a vertex algebra and let $\sigma$ be an automorphism of $V$ with a positive integer $T$ such that $\sigma^T=1$. 
A {\em right $\sigma$-twisted $V$-module} is a vector space $W$ equipped with a linear map
\begin{eqnarray*}
Y_{W}(\cdot,x):&&V\rightarrow (\End W)[[x^{\frac{1}{T}},x^{-\frac{1}{T}}]],\nonumber\\
&&v\mapsto Y_{W}(v,x),
\end{eqnarray*}
satisfying the conditions that $Y_{W}({\bf 1},x)=1$,
\begin{eqnarray}
Y_{W}(v,x)w\in W((x^{-\frac{1}{T}}))\   \   \  \mbox{ for }v\in V,\ w\in W,
\end{eqnarray}
and that 
\begin{eqnarray}\label{ehl}
& &x_{0}^{-1}\delta\left(\frac{x_{1}-x_{2}}{x_{0}}\right)
Y_{W}(v,x_{2})Y_{W}(u,x_{1})w
-x_{0}^{-1}\delta\left(\frac{x_{2}-x_{1}}{-x_{0}}\right)
Y_{W}(u,x_{1})Y_{W}(v,x_{2})w\  \   \nonumber\\
&&\hspace{2cm}=x_{2}^{-1}\delta\left(\frac{x_{1}-x_{0}}{x_{2}}\right)\left(\frac{x_{1}-x_{0}}{x_{2}}\right)^{\frac{j}{T}}
Y_{W}(Y(u,x_{0})v,x_{2})w
\end{eqnarray}
for $u\in V^{j},\ v\in V,\ w\in W$. }
\ed

The following is a reformulation of Proposition \ref{ptdual}:

\bl{lequiv-left-right-module}
Assume that $V$ is a conformal vertex algebra. Let $W$ be a vector space with a linear map 
$$Y_{W}(\cdot,x):\  V\rightarrow (\End W)[[x^{\frac{1}{T}},x^{-\frac{1}{T}}]].$$
For $v\in V$, set 
$$Y_{W}^{o}(v,x)=Y_{W}(e^{xL(1)}(-x^{-2})^{L(0)}v,x^{-1}).$$
 Then $(W,Y_{W})$ is a (left) $\sigma$-twisted $V$-module if and only if 
 $(W,Y^{o}_{W})$ is a right $\sigma$-twisted $V$-module with $V$ viewed as a vertex algebra.
\el

Assume that $(W,Y_{W})$ is a right $\sigma$-twisted module for a vertex algebra $V$. Note that for any $v\in V,\ z_0\in \C$,
$$Y_{W}(v,x+z_0)=e^{z_0\frac{d}{dx}}Y_W(v,x)\  \  \text{ exists in }\Hom (W,W((x^{-\frac{1}{T}}))).$$

The following is a twisted analogue of an observation in \cite{li-reg}:

\bl{lright-translation}
Let $V$ be a vertex algebra, let $(W,Y_{W})$ be a right $\sigma$-twisted $V$-module, 
and let $z_0$ be any complex number.  
For $v\in V$, set
\begin{eqnarray}
Y_{W}^{(z_0)}(v,x)=Y_{W}(v,x+z_0).
\end{eqnarray}
Then $(W,Y_{W}^{(z_0)})$ carries the structure of a right $\sigma$-twisted $V$-module.
\el

We formulate the following straightforward result to conclude this section:

\bl{lemma-twisted-module-tp}
Let $\sigma_1$ and $\sigma_2$ be finite-order automorphisms of vertex algebras $U_{1}$ and $U_{2}$, respectively. 
Set $\sigma=\sigma_1\otimes \sigma_2$, an automorphism of vertex algebra $U_1\otimes U_2$.
Suppose $W$ is a vector space with $\sigma_i$-twisted $U_{i}$-module structures $Y_W^{i}(\cdot,x)$
for $i=1,2$ such that
\begin{eqnarray}
Y_W^{1}(u_1,x_1)Y_W^2(u_2,x_2)=Y_W^2(u_2,x_2)Y_W^{1}(u_1,x_1)\    \    \   \mbox{ for }u_1\in U_1,\ u_2\in U_2. 
\end{eqnarray}
Then $W$ is a $\sigma$-twisted $U_1\otimes U_2$-module with $Y(u_1\otimes u_2,x)=Y_W^1(u_1,x)Y_W^2(u_2,x)$
for $u_1\in U_1,\ u_2\in U_2$.
\el

\section{$\sigma_1\times \sigma_{2}$-twisted $V\otimes V$-module 
${\mathfrak{D}}_{\sigma_1,\sigma_2}^{(z)}(W)$}

This section is the core of the paper, in which we present the construction of the $\sigma_1\times \sigma_{2}$-twisted 
$V\otimes V$-module ${\mathfrak{D}}_{\sigma_1,\sigma_2}^{(z)}(W)$.
Throughout this section, we assume that $V$ is a vertex algebra, $\sigma, \sigma_1,\sigma_2$ are finite-order automorphisms of $V$, 
and $T$ is a positive integer, such that
\begin{eqnarray}
\sigma=(\sigma_1\sigma_2)^{-1},\quad \sigma_1\sigma_2=\sigma_2\sigma_1,\quad \sigma_1^T=1=\sigma_2^T.
\end{eqnarray}
For any right $\sigma$-twisted $V$-module $(W,Y_W)$ and any nonzero complex number $z$, in this section
 we construct a $\sigma_1\otimes \sigma_2$-twisted $V\otimes V$-module 
${\mathfrak{D}}_{\sigma_1,\sigma_2}^{(z)}(W)$ inside $W^{*}$.

First of all, following \cite{hl1}, we fix a branch of logarithmic function with
\begin{eqnarray}
\log \xi=\ln |\xi|+i{\rm arg}\; \xi, \quad 0\le {\rm arg}\; \xi<2\pi 
\quad \text{for }\xi\in \mathbb{C}^{\times}.
\end{eqnarray}
Then define 
\begin{eqnarray}
\xi^{\alpha}=e^{\alpha\log \xi}
\end{eqnarray}
 for $\xi\in \mathbb{C}^{\times},\ \alpha\in \mathbb{C}$, in particular, 
\begin{eqnarray}
(-1)^{\alpha}=e^{\alpha \pi i}.
\end{eqnarray}
We have $\xi^{\alpha}\xi^{\beta}=\xi^{\alpha+\beta}$ for $\alpha,\beta\in \mathbb{C}$, though 
it is {\em not} true that $(\xi_1\xi_2)^{\alpha}=\xi_1^{\alpha}\xi_2^{\alpha}$ for $\xi_1,\xi_2\in \C^{\times},\ \alpha\in \C$.
We especially mention that $(-\xi)^{\alpha}\ne (-1)^{\alpha}\xi^{\alpha}$ for general $\xi\in \C^{\times},\ \alpha\in \C$.

For $\xi\in \C^{\times},\ \alpha\in \C,$ define
\begin{eqnarray}
& &(x\pm \xi)^{\alpha}=\sum_{j\ge 0}{\alpha\choose j}(\pm \xi)^{j}x^{\alpha-j}
\in x^{\alpha}{\C}[[x^{-1}]],\\
& &(\xi\pm x)^{\alpha}=\sum_{j\ge 0}{\alpha\choose j}\xi^{\alpha-j}(\pm x)^{j}
\in {\C}[[x]].
\end{eqnarray}
We have
\begin{eqnarray*}
&&(x\pm \xi)^{\alpha}(x\pm \xi)^{\beta}=(x\pm \xi)^{\alpha+\beta},\\
&&(\xi\pm x)^{\alpha}(\xi\pm x)^{\beta}=(\xi\pm x)^{\alpha+\beta}
\end{eqnarray*}
for $\alpha,\beta\in \C$. We shall often use the following relation
\begin{eqnarray}\label{edeltasub1}
x_{2}^{-1}\delta\left(\frac{x_{1}-x_{0}}{x_{2}}\right)(x_{2}-z)^{\alpha}
=x_{2}^{-1}\delta\left(\frac{x_{1}-x_{0}}{x_{2}}\right)
\left(\frac{x_{1}-x_{0}}{x_{2}}\right)^{-\alpha}(x_{1}-x_{0}-z)^{\alpha}.
\end{eqnarray}

For $j_1,j_2\in \Z$, set
\begin{eqnarray}
V^{(j_1,j_2)}=\{ v\in V\ |\ \sigma_r(v)=e^{2\pi i j_r/T},\ r=1,2\}.
\end{eqnarray}
Note that for $v\in V^{(j_1,j_2)}$, with $\sigma=(\sigma_1\sigma_2)^{-1}$ we have 
$$\sigma(v)=e^{-2\pi i (j_1+j_2)/T}v.$$

Let  $(W,Y_W)$ be a right $\sigma$-twisted $V$-module, which is fixed throughout this section.
 For $u\in V^{(j_1,j_2)},\ v\in V$ with $j_1,j_2\in \Z$, we have
\begin{eqnarray}\label{e-right-jacobiy}
& &x_{0}^{-1}\delta\left(\frac{x_{1}-x_{2}}{x_{0}}\right)
Y_W(v,x_{2})Y_W(u,x_{1})
-x_{0}^{-1}\delta\left(\frac{x_{2}-x_{1}}{-x_{0}}\right)
Y_W(u,x_{1})Y_W(v,x_{2})\nonumber\\
&&\hspace{2cm}   =x_{2}^{-1}\delta\left(\frac{x_{1}-x_{0}}{x_{2}}\right)\left(\frac{x_{1}-x_{0}}{x_{2}}\right)^{-\frac{j_1+j_2}{T}}
Y_W(Y(u,x_{0})v,x_{2}).
\end{eqnarray}
Note that this opposite twisted Jacobi identity implies
\begin{eqnarray}
x^{{j_1+j_2\over T}}Y_W(u,x)\in (\End W)[[x,x^{-1}]].
\end{eqnarray}

Define a linear map
$$Y_{W}^{*}(\cdot,x):\  V\rightarrow (\End\; W^*)[[x^{\frac{1}{T}},x^{-\frac{1}{T}}]]$$ 
by
\begin{eqnarray}
\< Y_{W}^{*}(v,x)f,w\>=\<f,Y_{W}(v,x)w\>\   \   \   \   \mbox{ for }v\in V,\ f\in W^*,\ w\in W.
\end{eqnarray}
For $u\in V^{(j_1,j_2)},\ v\in V,\ \alpha\in W^{*}$ with $j_1,j_2\in \Z$, we have
\begin{eqnarray}\label{ejacobiy*}
& &x_{0}^{-1}\delta\left(\frac{x_{1}-x_{2}}{x_{0}}\right)
Y_W^{*}(u,x_{1})Y_W^{*}(v,x_{2})\alpha
-x_{0}^{-1}\delta\left(\frac{x_{2}-x_{1}}{-x_{0}}\right)
Y_W^{*}(v,x_{2})Y_W^{*}(u,x_{1})\alpha\nonumber\\
&&\hspace{2cm}   =x_{2}^{-1}\delta\left(\frac{x_{1}-x_{0}}{x_{2}}\right)\left(\frac{x_{1}-x_{0}}{x_{2}}\right)^{-\frac{j_1+j_2}{T}}
Y_W^{*}(Y(u,x_{0})v,x_{2})\alpha.
\end{eqnarray}
This particularly implies the following commutator formula
\begin{eqnarray}\label{ecommutator*}
& &\quad Y_W^{*}(u,x_{1})Y_W^{*}(v,x_{2})\alpha-Y_W^{*}(v,x_{2})Y_W^{*}(u,x_{1})\alpha\nonumber\\
&&=\Res_{x_{0}}x_{2}^{-1}\delta\left(\frac{x_{1}-x_{0}}{x_{2}}\right)\left(\frac{x_{1}-x_{0}}{x_{2}}\right)^{-\frac{j_1+j_2}{T}}
Y_W^{*}(Y(u,x_{0})v,x_{2})\alpha.
\end{eqnarray}

\bd{dtwisted0} 
{\em Let $(W,Y_W)$ be a right $\sigma$-twisted $V$-module and let $z$ be a nonzero complex number. 
Denote by ${\mathfrak{D}}_{\sigma_1,\sigma_2}^{(z)}(W)$ the set of all linear 
functionals $\alpha$ on $W$, satisfying the condition that 
for every $v\in V^{(j_1,j_2)}$ with  $j_1,j_2\in \Z$, there exists $k\in {\N}$ such that
\begin{eqnarray}\label{etd}
x^{{j_2\over T}}(x-z)^{k+{j_1\over T}}Y_W^{*}(v,x)\alpha\in W^{*}((x)).
\end{eqnarray}}
\ed

From definition, ${\mathfrak{D}}_{\sigma_1,\sigma_2}^{(z)}(W)$ is a subspace of $W^{*}$.
Notice that if (\ref{etd}) holds for some nonnegative integer $k$, then it holds  
for any integer bigger than $k$.

The following follows immediately from the definition:

\bl{lbasic0}
Let $\alpha\in W^{*}$.
Then $\alpha\in {\mathfrak{D}}_{\sigma_1,\sigma_2}^{(z)}(W)$ if and only if for every $v\in V^{(j_1,j_2)}$ with $j_1,j_2\in \Z$, 
there exist $r,s\in {\N}$ such that
\begin{eqnarray}\label{char-property}
x^{r+{j_2\over T}}(x-z)^{s+{j_1\over T}}Y_W^{*}(v,x)\alpha\in W^{*}[[x]].
\end{eqnarray}
\el

Note that condition (\ref{char-property}) is equivalent to
\begin{eqnarray}
x^{r+{j_2\over T}}(x-z)^{s+{j_1\over T}}\< Y_W^{*}(v,x)\alpha,w\> \in {\C}[x] 
\end{eqnarray}
for all $w\in W$ as $x^{{1\over T}(j_1+j_2)}Y_W(v,x)w\in W((x^{-1}))$.

\bl{lbasic1} 
Let $ \alpha\in {\mathfrak{D}}_{\sigma_1,\sigma_2}^{(z)}(W),\  u\in V^{(j_{1},j_2)},\  v\in V^{(r_1,r_2)}$ with $j_1,j_2,r_1,r_2\in \Z$
and let $n\in \Z$. 
Then there exists $k\in {\N}$ such that 
\begin{eqnarray}
(x_{1}-z)^{k+{j_{1}\over T}}(x_{2}-z)^{n+{r_1\over T}}
Y_W^{*}(u,x_{1})Y_W^{*}(v,x_{2})\alpha\in x_2^{-\frac{r_2}{T}}x_{1}^{-\frac{j_2}{T}}W^{*}((x_{1}))[[x_{2},x_{2}^{-1}]].
\end{eqnarray}
Furthermore, we have 
\begin{eqnarray}
x^{r_2\over T}(x-z)^{n+{r_1\over T}}Y_W^{*}(v,x)\alpha
\in {\mathfrak{D}}_{\sigma_1,\sigma_2}^{(z)}(W)[[x,x^{-1}]].
\end{eqnarray}
\el

\begin{proof} Let $k$ be a nonnegative integer to be determined. Using the commutator formula (\ref{ecommutator*}) 
with $(u,x_1,x_0)\leftrightarrow (v,x_2,x_0)$, we have
\begin{eqnarray}
& &(x_{1}-z)^{k+{j_{1}\over T}}(x_2-z)^{n+{r_{1}\over T}}
Y_W^{*}(u,x_{1})Y_W^{*}(v,x_{2})\alpha\nonumber\\
&=&(x_{1}-z)^{k+{j_{1}\over T}}(x_2-z)^{n+{r_1\over T}}
Y_W^{*}(v,x_{2})Y_W^{*}(u,x_{1})\alpha\nonumber\\
& &-\Res_{x_{0}}x_{1}^{-1}\delta\left(\frac{x_{2}-x_{0}}{x_{1}}\right)\left(\frac{x_{2}-x_{0}}{x_{1}}\right)^{-\frac{r_1+r_2}{T}}
(x_{1}-z)^{k+{j_{1}\over T}}(x_2-z)^{n+{r_1\over T}}
Y_W^{*}(Y(v,x_{0})u,x_{1})\alpha \nonumber\\
&=&(x_{1}-z)^{k+{j_{1}\over T}}(x_2-z)^{n+{r_1\over T}}
Y_W^{*}(v,x_{2})Y_W^{*}(u,x_{1})\alpha\nonumber\\
& &-\Res_{x_{0}}x_{2}^{-1}\delta\left(\frac{x_{1}+x_{0}}{x_{2}}\right)\left(\frac{x_{1}+x_{0}}{x_{2}}\right)^{\frac{r_1+r_2}{T}}
(x_{1}-z)^{k+{j_{1}\over T}}(x_2-z)^{n+{r_1\over T}}
Y_W^{*}(Y(v,x_{0})u,x_{1})\alpha \nonumber\\
&=&(x_{1}-z)^{k+{j_{1}\over T}}(x_{2}-z)^{n+{r_1\over T}}
Y_W^{*}(v,x_{2})Y^{*}(u,x_{1})\alpha\nonumber\\
& &-\Res_{x_{0}}x_{2}^{-1}\delta\left(\frac{x_{1}+x_{0}}{x_{2}}\right)
(x_{1}-z)^{k+{j_{1}\over T}}(x_{1}+x_{0}-z)^{n+{r_1\over T}}
(x_{1}+x_{0})^{{r_2\over T}-n}x_2^{n-\frac{r_2}{T}}\nonumber\\
&& \quad \cdot 
Y_W^{*}(Y(v,x_{0})u,x_{1})\alpha \nonumber\\
&=&(x_{1}-z)^{k+{j_{1}\over T}}(x_{2}-z)^{n+{r_1\over T}}
Y_W^{*}(v,x_{2})Y^{*}(u,x_{1})\alpha\nonumber\\
& &-\sum_{r,s,t\in {\N}}\binom{n+\frac{r_1}{T}}{s}\binom{\frac{r_2}{T}-n}{t}
(x_{1}-z)^{k+n-s+{j_{1}+r_{1}\over T}}x_{1}^{{r_2\over T}-n-t}x_2^{n-\frac{r_2}{T}}Y_W^{*}(v_{r+s+t}u,x_{1})\alpha
 \nonumber\\
&&\   \   \   \cdot {1\over r!}\left({\partial\over\partial x_{1}}\right)^{r}
x_{2}^{-1}\delta\left({x_{1}\over x_{2}}\right).\label{e3.7}
\end{eqnarray}
Since $v_{m}u=0$ for all but finitely many $m\in {\N}$, we can choose
$k\in {\N}$ such that
$$x^{\frac{j_2}{T}}(x-z)^{k+{j_{1}\over T}}Y^{*}(u,x)\alpha,\quad
x^{\frac{j_2+r_2}{T}}(x-z)^{k+n-s+{j_{1}+r_1\over T}}Y_W^{*}(v_{m}u,x)\alpha\in W^{*}((x))$$
for all $m\in {\N}$. Then by (\ref{e3.7}) we obtain
$$x_1^{\frac{j_2}{T}} x_2^{\frac{r_2}{T}}(x_{1}-z)^{k+{j_{1}\over T}}(x_{2}-z)^{n+{r_1\over T}}
Y_W^{*}(u,x_{1})Y^{*}(v,x_{2})\alpha
\in W^{*}((x_{1}))[[x_{2},x_{2}^{-1}]],$$
as desired. The ``furthermore" assertion is clear.
\end{proof}

Let $\alpha\in {\mathfrak{D}}_{\sigma_1,\sigma_2}^{(z)}(W),\ v\in V^{(j_1,j_2)}$ with $j_1,j_2\in \Z$. 
Let $k\in {\N}$ be such that (\ref{etd}) holds. From Lemma \ref{lbasic1}, we have
\begin{eqnarray}
(x-z)^{k+{j_1\over T}}Y_W^{*}(v,x)\alpha\in x^{-{j_2\over T}}{\mathfrak{D}}_{\sigma_1,\sigma_2}^{(z)}(W)((x)),
\end{eqnarray}
so that
\begin{eqnarray}\label{ebrac}
(z-x)^{-k-{j_1\over T}}\left[(x-z)^{k+{j_1\over T}}Y_W^{*}(v,x)\alpha\right]
\    \    \mbox{exists in }x^{-{j_2\over T}}{\mathfrak{D}}_{\sigma_1,\sigma_2}^{(z)}(W)((x)).
\end{eqnarray}
Note that here {\em we must use the left and right brackets} because
$(z-x)^{-k-{j_1\over T}}Y_W^{*}(v,x)\alpha$ does not exist in $x^{-{j_2\over T}}W^{*}[[x,x^{-1}]]$ 
in general. Furthermore, for any $s\in {\N}$, 
we have
\begin{eqnarray}
& &e^{-(k+s+{j_1\over T})\pi i}(z-x)^{-k-s-{j_1\over T}}\left[(x-z)^{k+s+{j_1\over T}}Y_W^{*}(v,x)\alpha\right]
\nonumber\\
&=&e^{-(k+s+{j_1\over T})\pi i}(z-x)^{-k-s-{j_1\over T}}\left\{(x-z)^{s}\left[(x-z)^{k+{j_1\over T}}
Y_W^{*}(v,x)\alpha\right]\right\}
\nonumber\\
&=&e^{-(k+{j_1\over T})\pi i}(z-x)^{-k-{j_1\over T}}\left[(x-z)^{k+{j_1\over T}}Y_W^{*}(v,x)\alpha\right].
\end{eqnarray}
Therefore, the expression
\begin{eqnarray}
e^{-(k+{j_1\over T})\pi i}(z-x)^{-k-{j_1\over T}}\left[(x-z)^{k+{j_1\over T}}Y_W^{*}(v,x)\alpha\right]
\end{eqnarray}
does not depend on the choice of  $k$. On the other hand, it depends on $j_1$ modulo $T$. 

\bd{dtwisted1} 
{\em Let  $\alpha\in {\mathfrak{D}}^{(z)}_{\sigma_1,\sigma_2}(W)$. For $v\in V^{(j_1,j_2)}$ with $j_1,j_2\in \Z$, define 
\begin{eqnarray}
Y^{R}_{\sigma}(v,x)\alpha=e^{-(k+{j_1\over T})\pi i}
(z-x)^{-k-{j_1\over T}}\left((x-z)^{k+{j_1\over T}}Y_W^{*}(v,x)\alpha\right),
\end{eqnarray}
which lies in $x^{-{j_2\over T}}{\mathfrak{D}}^{(z)}_{\sigma_1,\sigma_2}(W)((x))$,
where $k$ is any nonnegative integer such that  (\ref{etd}) holds.}
\ed

As an  immediate consequence of  Definition \ref{dtwisted1}, we have:

\bc{cbasic} 
For $\alpha\in {\mathfrak{D}}_{\sigma_1,\sigma_2}^{(z)}(W),\ v\in V^{(j_1,j_2)}$ with $j_1,j_2\in \Z$,
\begin{eqnarray}\label{e3.11}
e^{(k+{j_1\over T})\pi i}(z-x)^{k+{j_1\over T}}Y_{\sigma}^{R}(v,x)\alpha=(x-z)^{k+{j_1\over T}}Y_W^{*}(v,x)\alpha
\end{eqnarray}
whenever $x^{j_2\over T}(x-z)^{k+{j_1\over T}}Y_W^{*}(v,x)\alpha\in W^{*}((x))$ with $k\in {\N}$.
Furthermore, we have
\begin{eqnarray}\label{e3.12}
e^{(k+{j_1\over T})\pi i}(z-x)^{k+{j_1\over T}}\<Y_{\sigma}^{R}(v,x)\alpha,w\>
=(x-z)^{k+{j_1\over T}}\<\alpha,Y_W(v,x)w\>
\end{eqnarray}
for all $w\in W$.
\ec

The following is the first key result:

\bt{thm-twistedr}
The pair $({\mathfrak{D}}^{(z)}_{\sigma_1,\sigma_2}(W),Y^{R}_{\sigma_2})$ carries the structure of a
$\sigma_{2}$-twisted $V$-module.
\et

\begin{proof} Let $u\in V^{(j_{1},j_2)},\  v\in V^{(r_1,r_2)},\  \alpha\in {\mathfrak{D}}^{(z)}_{\sigma_1,\sigma_2}(W)$ with $j_1,j_2,r_1,r_2\in \Z$.
There exists $k'\in {\N}$ such that
$$e^{(k'+{r_1\over T})\pi i}(z-x_{2})^{k'+{r_{1}\over T}}Y^{R}_{\sigma}(v,x_{2})\alpha
=(x_{2}-z)^{k'+{r_1\over T}}Y_W^{*}(v,x_{2})\alpha.$$
From Lemma \ref{lbasic1} we may choose $k'$ such that
\begin{eqnarray}
& &e^{({j_1\over T}+{r_1\over T})\pi i}(z-x_{1})^{k'+{j_{1}\over T}}(z-x_{2})^{k'+{r_1\over T}}
Y^{R}_{\sigma}(u,x_{1})Y^{R}_{\sigma}(v,x_{2})\alpha\nonumber\\
&=&e^{(k'+{j_{1}\over T})\pi i}(z-x_{1})^{k'+{j_{1}\over T}}(x_{2}-z)^{k'+{r_1\over T}}
Y^{R}_{\sigma}(u,x_{1})Y_W^{*}(v,x_{2})\alpha\nonumber\\
&=&(x_{1}-z)^{k'+{j_{1}\over T}}(x_{2}-z)^{k'+{r_1\over T}}
Y_W^{*}(u,x_{1})Y_W^{*}(v,x_{2})\alpha.
\end{eqnarray}
Clearly, this assertion with $(u,x_{1},v,x_{2})$ being replaced by $(v,x_{2},u,x_{1})$
also holds. Then we may choose $k'$ such that the following also holds:
\begin{eqnarray}
& &e^{({j_1\over T}+{r_1\over T})\pi i}(z-x_{1})^{k'+{j_{1}\over T}}(z-x_{2})^{k'+{r_1\over T}}
Y^{R}_{\sigma}(v,x_{2})Y^{R}_{\sigma}(u,x_{1})\alpha\nonumber\\
&=&(x_{1}-z)^{k'+{j_{1}\over T}}(x_{2}-z)^{k'+{r_1\over T}}
Y_W^{*}(v,x_{2})Y_W^{*}(u,x_{1})\alpha.
\end{eqnarray}
Using these relations we get
\begin{eqnarray}\label{estep1}
& &x_{0}^{-1}\delta\left(\frac{x_{1}-x_{2}}{x_{0}}\right)e^{({j_1\over T}+{r_1\over T})\pi i}
(z-x_{1})^{k'+{j_{1}\over T}}(z-x_{2})^{k'+{r_1\over T}}
Y^{R}_{\sigma}(u,x_{1})Y^{R}_{\sigma}(v,x_{2})\alpha\nonumber\\
& &\  \  \  \   -x_{0}^{-1}\delta\left(\frac{x_{2}-x_{1}}{-x_{0}}\right)e^{({j_1\over T}+{r_1\over T})\pi i}
(z-x_{1})^{k'+{j_{1}\over T}}(z-x_{2})^{k'+{r_1\over T}}
Y^{R}_{\sigma}(v,x_{2})Y^{R}_{\sigma}(u,x_{1})\alpha\nonumber\\
&=&x_{0}^{-1}\delta\left(\frac{x_{1}-x_{2}}{x_{0}}\right)
(x_{1}-z)^{k'+{j_{1}\over T}}(x_{2}-z)^{k'+{r_1\over T}}
Y_W^{*}(u,x_{1})Y_W^{*}(v,x_{2})\alpha\nonumber\\
& &-x_{0}^{-1}\delta\left(\frac{x_{2}-x_{1}}{-x_{0}}\right)
(x_{1}-z)^{k'+{j_{1}\over T}}(x_{2}-z)^{k'+{r_1\over T}}
Y_W^{*}(v,x_{2})Y_W^{*}(u,x_{1})\alpha\nonumber\\
&=&x_{2}^{-1}\delta\left(\frac{x_{1}-x_{0}}{x_{2}}\right)\left(\frac{x_{1}-x_{0}}{x_{2}}\right)^{-\frac{j_1+j_2}{T}}
(x_{1}-z)^{k'+{j_{1}\over T}}(x_{2}-z)^{k'+{r_1\over T}}
Y_W^{*}(Y(u,x_{0})v,x_{2})\alpha\nonumber\\
&=&x_{2}^{-1}\delta\left(\frac{x_{1}-x_{0}}{x_{2}}\right)
\left(\frac{x_{1}-x_{0}}{x_{2}}\right)^{-\frac{j_1+j_2}{T}-k'-{r_1\over T}}
(x_{1}-z)^{k'+{j_{1}\over T}}(x_{1}-x_{0}-z)^{k'+{r_1\over T}}\nonumber\\
&&\hspace{1cm}\cdot Y_W^{*}(Y(u,x_{0})v,x_{2})\alpha\nonumber\\
&=&x_{2}^{-1}\delta\left(\frac{x_{1}-x_{0}}{x_{2}}\right)
\left(\frac{x_{1}-x_{0}}{x_{2}}\right)^{-{j_1+j_2+r_1\over T}}
(x_{1}-z)^{k'+{j_{1}\over T}}(x_{1}-x_{0}-z)^{k'+{r_1\over T}}\nonumber\\
&&\hspace{1cm}\cdot Y_W^{*}(Y(u,x_{0})v,x_{2})\alpha,
\end{eqnarray}
where we use (\ref{edeltasub1}) for the last two equalities.

Let $m$ be {\em an arbitrarily fixed} integer. Notice that $u_{n}v\in V^{(j_{1}+r_1,j_{2}+r_2)}$ for $n\in \Z$ and
there are finitely many $n\ge m$
such that $u_{n}v\ne 0$. Then there exists $k''\in {\N}$ ({\em depending} on $m$) such that
$$(x-z)^{k''+{j_{1}+r_{1}\over T}}Y_W^{*}(u_{n}v,x)\alpha
\in x^{-{j_{2}+r_{2}\over T}}W^{*}((x))$$
 and hence
\begin{eqnarray*}
e^{(k''+{j_{1}+r_{1}\over T})\pi i} (z-x)^{k''+{j_{1}+r_{1}\over T}}Y_{\sigma}^{R}(u_{n}v,x)\alpha=
(x-z)^{k''+{j_{1}+r_{1}\over T}}Y_W^{*}(u_{n}v,x)\alpha
\end{eqnarray*}
for {\em all }$n\ge m$. Thus
\begin{eqnarray}\label{eadded-new}
&&\Res_{x_0}x_0^{m+s}e^{(k''+{j_{1}+r_1\over T})\pi i}(z-x)^{k''+{j_{1}+r_1\over T}}Y_{\sigma}^{R}(Y(u,x_0)v,x)\alpha\nonumber\\
&=&\Res_{x_0}x_0^{m+s}(x-z)^{k''+{j_{1}+r_1\over T}}Y_W^{*}(Y(u,x_0)v,x)\alpha
\end{eqnarray}
for all $s\ge 0$.
Let $r\in {\N}$ be such that $u_{n}v=0$ for $n\ge r$.
By choosing $k\in {\N}$ such that $k\ge k', k'', r$ and using (\ref{eadded-new}), we obtain
\begin{eqnarray}\label{estep2}
& &\Res_{x_{0}}x_{0}^{m}x_{2}^{-1}\delta\left(\frac{x_{1}-x_{0}}{x_{2}}\right)
\left(\frac{x_{1}-x_{0}}{x_{2}}\right)^{-{j_{1}+j_2+r_1\over T}}
(x_1-z)^{k+{j_{1}\over T}}(x_1-x_{0}-z)^{k+{r_1\over T}}\nonumber\\
&&\hspace{0.5cm}\cdot Y_W^{*}(Y(u,x_{0})v,x_{2})\alpha\nonumber\\
&=&\Res_{x_{0}}x_{0}^{m}x_{1}^{-1}\delta\left(\frac{x_2+x_{0}}{x_{1}}\right)
\left(\frac{x_2+x_{0}}{x_{1}}\right)^{{j_{1}+j_2+r_1\over T}-2k-{j_1\over T}-{r_1\over T}}
(x_2+x_0-z)^{k+{j_{1}\over T}}(x_2-z)^{k+{r_1\over T}}\nonumber\\
&&\hspace{0.5cm}\cdot Y_W^{*}(Y(u,x_{0})v,x_{2})\alpha\nonumber\\
&=&\Res_{x_{0}}\sum_{s\ge 0}{k+{j_1\over T}\choose s}x_{0}^{m+s}
(x_2-z)^{2k+{j_{1}+r_1\over T}-s}\cdot \nonumber\\
& &\cdot x_1^{-1}\delta\left(\frac{x_2+x_{0}}{x_1}\right)
\left(\frac{x_2+x_{0}}{x_1}\right)^{{j_2\over T}}
Y_W^{*}(Y(u,x_{0})v,x_{2})\alpha\nonumber\\
&=&\Res_{x_{0}}\sum_{s\ge 0}{k+{j_1\over T}\choose s}x_{0}^{m+s}
e^{(2k+{j_{1}+r_1\over T}-s)\pi i}  (z-x_2)^{2k+{j_{1}+r_1\over T}-s}\cdot \nonumber\\
& &\cdot x_1^{-1}\delta\left(\frac{x_2+x_{0}}{x_1}\right)
\left(\frac{x_2+x_{0}}{x_1}\right)^{{j_2\over T}}
Y_{\sigma}^{R}(Y(u,x_{0})v,x_{2})\alpha\nonumber\\
&=&\Res_{x_{0}}x_{0}^{m}e^{({j_1\over T}+{r_1\over T})\pi i}
(z-x_{2}-x_{0})^{k+{j_1\over T}}(z-x_{2})^{k+{r_{1}\over T}}\cdot \nonumber\\
& &\cdot x_1^{-1}\delta\left(\frac{x_2+x_{0}}{x_1}\right)
\left(\frac{x_2+x_{0}}{x_1}\right)^{{j_2\over T}}
Y_{\sigma}^{R}(Y(u,x_{0})v,x_{2})\alpha\nonumber\\
&=&\Res_{x_{0}}x_{0}^{m}e^{({j_1\over T}+{r_1\over T})\pi i}
(z-x_{1})^{k+{j_{1}\over T}}(z-x_{2})^{k+{r_1\over T}}
x_1^{-1}\delta\left(\frac{x_2+x_{0}}{x_1}\right)
\left(\frac{x_2+x_{0}}{x_1}\right)^{{j_2\over T}}\nonumber\\
&&\hspace{0.5cm}\cdot Y_{\sigma}^{R}(Y(u,x_{0})v,x_{2})\alpha.
\end{eqnarray}
Combining (\ref{estep1}) (with $k'$ being replaced by $k$) with (\ref{estep2}) 
we get
\begin{eqnarray}\label{eresjacobi}
& &\Res_{x_{0}}x_{0}^{m}x_{0}^{-1}\delta\left(\frac{x_{1}-x_{2}}{x_{0}}\right)
(z-x_{1})^{k+{j_{1}\over T}}(z-x_{2})^{k+{r_1\over T}}
Y^{R}_{\sigma}(u,x_{1})Y^{R}_{\sigma}(v,x_{2})\alpha\nonumber\\
& &-\Res_{x_{0}}x_{0}^{m}x_{0}^{-1}\delta\left(\frac{x_{2}-x_{1}}{-x_{0}}\right)
(z-x_{1})^{k+{j_{1}\over T}}(z-x_{2})^{k+{r_1\over T}}
Y^{R}_{\sigma}(v,x_{2})Y^{R}_{\sigma}(u,x_{1})\alpha\nonumber\\
&=&\Res_{x_{0}}x_{0}^{m}
(z-x_{1})^{k+{j_{1}\over T}}(z-x_{2})^{k+{r_1\over T}}
x_{1}^{-1}\delta\left(\frac{x_2+x_{0}}{x_1}\right)
\left(\frac{x_2+x_{0}}{x_1}\right)^{{j_2\over T}}\nonumber\\
&&\hspace{1cm}\cdot Y_{\sigma}^{R}(Y(u,x_{0})v,x_{2})\alpha.
\end{eqnarray}
Since $Y_{\sigma}^{R}(u,x_{1})$ and $Y_{\sigma}^{R}(v,x_{2})$ acting on ${\mathfrak{D}}^{(z)}_{\sigma_1,\sigma_2}(W)$
satisfy the truncation condition, we may multiply the both sides of (\ref{eresjacobi}) by 
$(z-x_{1})^{-k-{j_{1}\over T}}(z-x_{2})^{-k-{r_1\over T}}$
to get
\begin{eqnarray}\label{eresjacobi2}
& &\Res_{x_{0}}x_{0}^{m}x_{0}^{-1}\delta\left(\frac{x_{1}-x_{2}}{x_{0}}\right)
Y^{R}_{\sigma}(u,x_{1})Y^{R}_{\sigma}(v,x_{2})\alpha\nonumber\\
& &\   \   \   \    -\Res_{x_{0}}x_{0}^{m}x_{0}^{-1}\delta\left(\frac{x_{2}-x_{1}}{-x_{0}}\right)
Y^{R}_{\sigma}(v,x_{2})Y^{R}_{\sigma}(u,x_{1})\alpha\nonumber\\
&=&\Res_{x_{0}}x_{0}^{m}
x_1^{-1}\delta\left(\frac{x_2+x_{0}}{x_1}\right)
\left(\frac{x_2+x_{0}}{x_1}\right)^{{j_2\over T}}
Y_{\sigma}^{R}(Y(u,x_{0})v,x_{2})\alpha.
\end{eqnarray}
Now, in (\ref{eresjacobi2}) no terms except $x_{0}^{m}$ depend on $m$.
Since $m$ can be any integer,
dropping $\Res_{x_{0}}x_{0}^{m}$ off from (\ref{eresjacobi2}) we get the desired
$\sigma_{2}$-twisted Jacobi identity.
This completes the proof. 
\end{proof}

Note that for any $v\in V,\  w\in W$, $Y_{W}(v,x+z)w$ exists in $W((x^{-\frac{1}{T}}))$, hence
$Y_W^{*}(v,x+z)\alpha$ exists in $W^{*}[[x^{\frac{1}{T}},x^{-\frac{1}{T}}]]$ for any $v\in V,\ \alpha\in W^{*}$. Furthermore,
 for any real number $\mu$, $(x+z)^{\mu}Y_W^{*}(v,x+z)\alpha$ exists in $x^{\mu}W^{*}[[x^{\frac{1}{T}},x^{-\frac{1}{T}}]]$.
The following is another characterization of  ${\mathfrak{D}}^{(z)}_{\sigma_1,\sigma_2}(W)$:

\bl{pleftrightcondition}
Let $\alpha\in W^{*}$. Then $\alpha\in {\mathfrak{D}}^{(z)}_{\sigma_1,\sigma_2}(W)$ if and only if
for every $v\in V^{(j_1,j_2)}$ with $j_1,j_2\in \Z$, there exist nonnegative integers $r$ and $s$ such that
\begin{eqnarray}\label{eleftcon}
x^{s+{j_1\over T}}(x+z)^{r+{j_2\over T}}Y_W^{*}(v,x+z)\alpha\in W^{*}[[x]].
\end{eqnarray}
\el

\begin{proof} From Lemma \ref{lbasic0}, $\alpha\in {\mathfrak{D}}_{\sigma_1,\sigma_2}^{(z)}(W)$ if and only if
for every $v\in V^{(j_1,j_2)}$ with $j_1,j_2\in \Z$, there exist nonnegative integers $r,s$ such that
\begin{eqnarray}\label{epr1}
\<x^{r+{j_2\over T}}(x-z)^{s+{j_1\over T}}Y_W^{*}(v,x)\alpha,w\>\in {\C}[x]
\end{eqnarray}
for all $w\in W$. Similarly, 
(\ref{eleftcon}) holds if and only if there exist
$k_{1}, k_{2}\in {\N}$ such that
\begin{eqnarray}\label{epr2}
\<x^{k_{1}+{j_1\over T}}(x+z)^{k_{2}+{j_2\over T}}Y_W^{*}(v,x+z)\alpha,w\>\in {\C}[x]
\end{eqnarray}
for all $w\in W$. Applying substitution $x=x_{0}+z$ to (\ref{epr1}) we get  (\ref{epr2}) with $k_1=s$ and $k_2=r$.
On the other hand, applying substitution $x=x_{0}-z$ to (\ref{epr2}) we get (\ref{epr1}) with $s=k_1$ and $r=k_2$.
This completes the proof.
\end{proof}

We have the following analogue of Lemma \ref{lbasic1}: 

\bp{pyl0}
Let $\alpha\in {\mathfrak{D}}^{(z)}_{\sigma_1,\sigma_2}(W),\ v\in V^{(j_1,j_2)}$ with $j_1,j_2\in \Z$. Then
\begin{eqnarray}\label{enva}
x^{{j_1\over T}}(x+z)^{n+{j_2\over T}}Y_W^{*}(v,x+z)\alpha
\in {\mathfrak{D}}_{\sigma_1,\sigma_2}^{(z)}(W)[[x,x^{-1}]]
\end{eqnarray}
for all $n\in {\Z}$. 
\ep

\begin{proof} Let $u\in V^{(j_1,j_2)},\ v\in V^{(r_1,r_2)}$, $j_1,j_2,r_1,r_2\in \Z$. 
Applying ${\Res}_{x_{0}}$ to (\ref{ejacobiy*}) with $(u,v)$ being permuted we get
\begin{eqnarray*}
& &Y_W^{*}(v,x_{1}+z)Y_W^{*}(u,x_{2})\alpha
-Y_W^{*}(u,x_{2})Y_W^{*}(v,x_{1}+z)\alpha\nonumber\\
&=&\Res _{x_{0}}x_{2}^{-1}\delta\left(\frac{x_{1}+z-x_{0}}{x_{2}}\right)\left(\frac{x_{1}+z-x_{0}}{x_{2}}\right)^{-\frac{r_1+r_2}{T}}
Y_W^{*}(Y(v,x_{0})u,x_{2})\alpha\nonumber\\
&=&\Res _{x_{0}}x_{1}^{-1}\delta\left(\frac{x_{2}-z+x_{0}}{x_1}\right)\left(\frac{x_2-z+x_{0}}{x_1}\right)^{\frac{r_1+r_2}{T}}
Y_W^{*}(Y(v,x_{0})u,x_{2})\alpha.
\end{eqnarray*}
Let $m$ be an arbitrarily fixed integer. Then we have
\begin{eqnarray}\label{en*a2}
& &Y_W^{*}(u,x_{2})\Res_{x_{1}}x_{1}^{m+{r_1\over T}}(x_{1}+z)^{n+{r_2\over T}}
Y_W^{*}(v,x_{1}+z)\alpha
\nonumber\\
&=&\Res_{x_{1}}x_{1}^{m+{r_1\over T}}
(x_{1}+z)^{n+{r_2\over T}}Y_W^{*}(v,x_{1}+z)Y_W^{*}(u,x_{2})\alpha
\nonumber\\
& &-\Res_{x_{1}}\Res_{x_{0}}
x_1^{-1}\delta\left(\frac{x_2-z+x_{0}}{x_1}\right)\left(\frac{x_2-z+x_{0}}{x_1}\right)^{\frac{r_1+r_2}{T}}\nonumber\\
&&\quad \cdot x_{1}^{m+{r_1\over T}}
(x_{1}+z)^{n+{r_2\over T}}Y_W^{*}(Y(v,x_{0})u,x_{2})\alpha\nonumber\\
&=&\Res_{x_{1}}x_{1}^{m+{r_1\over T}}
(x_{1}+z)^{n+{r_2\over T}}Y_W^{*}(v,x_{1}+z)Y_W^{*}(u,x_{2})\alpha
\nonumber\\
& &-\Res_{x_{1}}\Res_{x_{0}}
x_{1}^{-1}\delta\left(\frac{x_{2}-z+x_{0}}{x_{1}}\right)(x_{2}-z+x_{0})^{\frac{r_1+r_2}{T}}
\left[x_{1}^{m-{r_2\over T}}
(x_{1}+z)^{n+{j_2\over T}}\right]\nonumber\\
&&\quad \cdot Y_W^{*}(Y(v,x_{0})u,x_{2})\alpha\nonumber\\
&=&\Res_{x_{1}}x_{1}^{m+{r_1\over T}}
(x_{1}+z)^{n+{r_2\over T}}Y_W^{*}(v,x_{1}+z)Y_W^{*}(u,x_{2})\alpha
\nonumber\\
& &-\Res_{x_{0}}\left[(x_{2}-z+x_{0})^{m+{r_1\over T}}
(x_{2}+x_{0})^{n+{r_2\over T}}\right]Y_W^{*}(Y(v,x_{0})u,x_{2})\alpha\nonumber\\
&=&\Res_{x_{1}}x_{1}^{m+{r_1\over T}}
(x_{1}+z)^{n+{r_2\over T}}Y_W^{*}(v,x_{1}+z)Y_W^{*}(u,x_{2})\alpha
\nonumber\\
& &-\sum_{0\le s, t\le r}{m+{r_1\over T}\choose s}{n+{r_2\over T}\choose t}
(x_{2}-z)^{m-s+{r_1\over T}}x_{2}^{n-t+{r_2\over T}}
Y_W^{*}(v_{s+t}u,x_{2})\alpha,
\end{eqnarray}
where $r$ is a nonnegative integer such that $v_{p}u=0$ for $p>r$.
Notice that in the previous equation, $x_{1}^{m-{r_2\over T}}(x_{1}+z)^{n+{r_2\over T}}$
involves only integer powers of $x_{1}$, so that we could use the 
delta function substitution.
Also note that $v_{p}u\in V^{(j_1+r_1,j_2+r_2)}$. 
Let $k\in {\N}$ be such that
$$(x_{2}-z)^{k+{j_{1}\over T}}Y_W^{*}(u,x_{2})\alpha\in x_2^{-\frac{j_2}{T}}W^{*}((x_{2}))$$
and such that $k+m-r\ge 0$ and
$$(x_{2}-z)^{k+m-r+{j_{1}+r_1\over T}}Y_W^{*}(v_{p}u,x_{2})\alpha\in x_2^{-\frac{j_2+r_2}{T}}W^{*}((x_{2}))$$
for all $0\le p\le r$. Then from (\ref{en*a2}) we obtain
$$(x_{2}-z)^{k+{j_{1}\over T}}Y_W^{*}(u,x_{2})
\Res_{x_{1}}x_{1}^{m+{r_1\over T}}(x_{1}+z)^{n+{r_2\over T}}Y_W^{*}(v,x_{1}+z)\alpha
\in x_2^{-\frac{j_2}{T}}W^{*}((x_{2})).$$
That is, 
$$\Res_{x}x^{m+{r_1\over T}}(x+z)^{n+{r_2\over T}}Y_W^{*}(v,x+z)\alpha
\in {\mathfrak{D}}_{\sigma_1,\sigma_2}^{(z)}(W).$$
Since $m$ is arbitrary, (\ref{enva}) follows. 
\end{proof}

{\bf A subtlety:} Note that the coefficient of $x^{q-\frac{j}{T}}$ in $(x+z)^{n-{j\over T}}Y_W^{*}(v,x+z)\alpha$ 
for $q\in \Z$ is an infinite but well-defined sum of $v_{m}^{*}\alpha$ with $m\in \Z$.

Now, we define another vertex operator map.

\bd{dleftaction}
{\em Define  a linear map
$$Y_{\sigma}^{L}(\cdot,x):V\rightarrow  (\End\; {\mathfrak{D}}^{(z)}_{\sigma_1,\sigma_2}(W))[[x,x^{-1}]]$$ 
by
\begin{eqnarray}
Y^{L}_{\sigma}(v,x)\alpha=(z+x)^{-l-{j_2\over T}}\left[ (x+z)^{l+{j_2\over T}}
Y_W^{*}(v,x+z)\alpha\right]
\end{eqnarray}
for $\alpha\in {\mathfrak{D}}^{(z)}_{\sigma_1,\sigma_2}(W), v\in V^{(j_1,j_2)}$, $j_1,j_2\in \Z$, 
where $l$ is any nonnegative integer such that
$$x^{{j_1\over T}}(x+z)^{l+{j_2\over T}}Y_W^{*}(v,x+z)\alpha \in W^{*}((x))$$
(cf. Proposition \ref{pyl0}).}
\ed

Note that the following identity holds (see \cite{hl1}):
\begin{eqnarray}\label{e3delta}
x^{-1}\delta\left(\frac{z+x_{0}}{x}\right)
=x_{0}^{-1}\delta\left(\frac{x-z}{x_{0}}\right)
-x_{0}^{-1}\delta\left(\frac{z-x}{-x_{0}}\right)
\end{eqnarray}
(with $z$ a nonzero complex number and with $x,x_0$ formal variables).
We also have the following delta-function substitution:
\begin{eqnarray}\label{esub}
& &x_{0}^{-1}\delta\left(\frac{z-x}{x_{0}}\right)
f(x_{0},x)=x_{0}^{-1}\delta\left(\frac{z-x}{x_{0}}\right)
f(z-x,x)
\end{eqnarray}
for $f(x_{0},x)\in \Hom (U_{1},U_{2}[x_{0},x_{0}^{-1},x,x^{-1}])$,
where $U_{1}$ and $U_{2}$ are any vector spaces. 

\bp{ptwistedleft}
The pair $({\mathfrak{D}}_{\sigma_1,\sigma_2}^{(z)}(W), Y_{\sigma}^{L})$ carries the structure of a 
$\sigma_1$-twisted $V$-module.
\ep

\begin{proof} For clarity, we locally denote ${\mathfrak{D}}_{\sigma_1,\sigma_2}^{(z)}(W)$ by ${\mathfrak{D}}_{\sigma_1,\sigma_2}^{(z)}(W,Y_{W})$,
and denote $Y_{\sigma}^R$, $Y_{\sigma}^L$ by $Y_{\sigma,z}^R$, $Y_{\sigma,z}^L$.
From Lemma \ref{lright-translation},  $(W,Y_{W}^{(z)})$ is a $\sigma$-twisted right $V$-module. 
Then by Theorem \ref{thm-twistedr} we have a
$\sigma_1$-twisted $V$-module $({\mathfrak{D}}_{\sigma_2,\sigma_1}^{(-z)}(W,Y_{W}^{(z)}),Y^R)$.
From Lemma \ref{pleftrightcondition}, we see that
\begin{eqnarray}
{\mathfrak{D}}_{\sigma_1,\sigma_2}^{(z)}(W,Y_{W})={\mathfrak{D}}_{\sigma_2,\sigma_1}^{(-z)}(W,Y_{W}^{(z)})
\end{eqnarray}
as vector spaces. Since $0\le \arg \xi<2\pi$ for $\xi\in \C^{\times}$, we have
$$\arg (-z)=\begin{cases}\arg z+ \pi\  \  \text{ if }\arg z<\pi\\
\arg z-\pi\  \  \text{ if }\arg z\ge \pi.\end{cases}$$ 
Then
$$(-z)^{\beta}=e^{\beta\pi i}z^{\beta}e^{\epsilon_z \beta 2\pi i}\quad \text{ for } \beta\in \C,$$
where $\epsilon_z=0,-1$ for $\arg z<\pi$ and $\arg z\ge \pi$, respectively.
  For $l,j\in \Z$, we have
$$e^{-(l+{j\over T})\pi i}(-z-x)^{-l-{j\over T}}=e^{-(1+\epsilon_z)\frac{j}{T}2\pi i}(z+x)^{-l-{j\over T}}. $$
If $\epsilon_z=-1$,  from Definitions \ref{dleftaction} and \ref{dtwisted1}
we see that $Y^{L}$ on ${\mathfrak{D}}_{\sigma_1,\sigma_2}^{(z)}(W,Y_{W})$ coincides with $Y^{R}$ on 
${\mathfrak{D}}_{\sigma_2,\sigma_1}^{(-z)}(W,Y_{W}^{(z)})$. Then by Theorem \ref{thm-twistedr}
$({\mathfrak{D}}_{\sigma_1,\sigma_2}^{(z)}(W), Y^{L})$ carries the structure of a 
$\sigma_1$-twisted $V$-module. Assume $\epsilon_z=0$.  From Definitions \ref{dleftaction} and \ref{dtwisted1} again,
we see that  for $v\in V^{(j_1,j_2)}$,
$$Y_{z}^L(v,x)=e^{\frac{j_2}{T}2\pi i}Y_{-z}^{R}(v,x)  =Y_{-z}^{R}(\sigma_2(v),x).$$
It then follows from Theorem \ref{thm-twistedr} and Lemma \ref{lem-sigma-tau}
that  $({\mathfrak{D}}_{\sigma_1,\sigma_2}^{(z)}(W), Y_{z}^{L})$ 
carries the structure of a $\sigma_1$-twisted $V$-module.
\end{proof}

Next, we are going to show that the two actions associated to $Y_{\sigma}^{R}$ and $Y_{\sigma}^{L}$ of $V$ 
on ${\mathfrak{D}}_{\sigma_1,\sigma_2}^{(z)}(W)$ commute.
First, we give a connection among $Y_W^{*}, Y_{\sigma}^{R}$, and $Y_{\sigma}^{L}$.

\bp{prelation}
Let $\alpha\in {\mathfrak{D}}_{\sigma_1,\sigma_2}^{(z)}(W),\ v\in V^{(j_1,j_2)}$ with $j_1,j_2\in \Z$. Then
\begin{eqnarray}\label{sjac}
& &x_{0}^{-1}\delta\left(\frac{x-z}{x_{0}}\right)\left(\frac{x-z}{x_{0}}\right)^{j_1\over T}
Y_W^{*}(v,x)\alpha  -e^{\frac{j_1}{T}\pi i}x_{0}^{-1}\delta\left(\frac{z-x}{-x_{0}}\right)
\left(\frac{z-x}{x_{0}}\right)^{j_1\over T}Y_{\sigma}^{R}(v,x)\alpha\nonumber\\
&&\hspace{1cm}=x^{-1}\delta\left(\frac{z+x_{0}}{x}\right)\left(\frac{z+x_{0}}{x}\right)^{{j_2\over T}}
Y_{\sigma}^{L}(v,x_{0})\alpha.
\end{eqnarray}
\ep

Before presenting the proof, we justify the existence of the three main terms in (\ref{sjac}).
Since $Y_{\sigma}^{R}(v,x)\alpha\in x^{-{j_2\over T}}{\mathfrak{D}}_{\sigma_1,\sigma_2}^{(z)}(W)((x))$, 
 $(z-x)^{n+{j_1\over T}}Y_{\sigma}^{R}(v,x)\alpha$ exists 
in $x^{-{j_2\over T}}{\mathfrak{D}}_{\sigma_1,\sigma_2}^{(z)}(W)((x))$ for every $n\in {\Z}$, so that
\begin{eqnarray*}
x_{0}^{-1}\delta\left(\frac{z-x}{-x_{0}}\right)\left(\frac{z-x}{x_{0}}\right)^{j_1\over T}
Y_{\sigma}^{R}(v,x)\alpha
\end{eqnarray*}
exists in $x_{0}^{-{j_1\over T}}x^{-{j_2\over T}}{\mathfrak{D}}_{\sigma_1,\sigma_2}^{(z)}(W)[[x_{0},x_{0}^{-1},x,x^{-1}]]$.
Similarly, the expression
\begin{eqnarray*}
x^{-1}\delta\left(\frac{z+x_{0}}{x}\right)\left(\frac{z+x_{0}}{x}\right)^{{j_2\over T}}
Y_{\sigma}^{L}(v,x_{0})\alpha
\end{eqnarray*}
 exists in $x_{0}^{-{j_1\over T}}x^{-{j_2\over T}}{\mathfrak{D}}_{\sigma_1,\sigma_2}^{(z)}(W)[[x_{0},x_{0}^{-1},x,x^{-1}]]$.
By Lemma \ref{lbasic1}, 
$(x-z)^{n+{j_1\over T}}Y_W^{*}(v,x)\alpha$ is a well-defined element of
$x^{-{j_2\over T}}{\mathfrak{D}}_{\sigma_1,\sigma_2}^{(z)}(W)[[x,x^{-1}]]$ for every $n\in {\Z}$, so that
\begin{eqnarray*}
x_{0}^{-1}\delta\left(\frac{x-z}{x_{0}}\right)
\left(\frac{x-z}{x_{0}}\right)^{j_1\over T}Y_W^{*}(v,x)\alpha
\end{eqnarray*}
is a well-defined element of 
$x_{0}^{-{j_1\over T}}x^{-{j_2\over T}}{\mathfrak{D}}_{\sigma_1,\sigma_2}^{(z)}(W)[[x_{0},x_{0}^{-1},x,x^{-1}]]$.

\begin{proof}   {\bf Proof of Proposition \ref{prelation}:}
First, by the definitions of $Y_{\sigma}^{R}$ and $Y_{\sigma}^{L}$, we have
\begin{eqnarray}
& &(x-z)^{k+{j_1\over T}}Y_W^{*}(v,x)\alpha=e^{(k+{j_1\over T})\pi i}(z-x)^{k+{j_1\over T}}Y_{\sigma}^{R}(v,x)\alpha,
\label{ecommr}\\
& &(x_{0}+z)^{k+{j_2\over T}}Y_W^{*}(v,x_{0}+z)\alpha
=(z+x_{0})^{k+{j_2\over T}}Y_{\sigma}^{L}(v,x_{0})\alpha\label{eassocr}
\end{eqnarray}
for some $k\in {\N}$. These are similar to
the weak commutativity relation (\ref{etwistedweakcomm}) and the weak associativity 
relation (\ref{etwistedweakassoc}).
In the following, we prove that (\ref{sjac}) follows from (\ref{ecommr}) and (\ref{eassocr})
just as the Jacobi identity follows from the weak commutativity and 
associativity in three formal variables (cf. \cite{li-local}).

In view of Lemma \ref{lbasic0}, in the above we may choose $k$ so large  that
\begin{eqnarray}\label{enon}
x^{k+{j_2\over T}}(x-z)^{k+{j_1\over T}}Y_W^{*}(v,x)\alpha\in \Hom (W, {\C}[x]).
\end{eqnarray}
Then using delta function properties (\ref{e3delta}) and (\ref{esub}) we obtain
\begin{eqnarray*}
& &x_{0}^{-1}\delta\left(\frac{x-z}{x_{0}}\right)\left(\frac{x-z}{x_{0}}\right)^{j_1\over T}
x^{k+{j_2\over T}}x_{0}^{k}Y_W^{*}(v,x)\alpha\nonumber\\
& &\   \   \   \   -e^{\frac{j_1}{T}\pi i}x_{0}^{-1}\delta\left(\frac{z-x}{-x_{0}}\right)\left(\frac{z-x}{x_{0}}\right)^{j_1\over T}
x^{k+{j_2\over T}}x_{0}^{k}Y_{\sigma}^{R}(v,x)\alpha\nonumber\\
&=&x_{0}^{-1}\delta\left(\frac{x-z}{x_{0}}\right)\left(\frac{x-z}{x_{0}}\right)^{j_1\over T}
x^{k+{j_2\over T}}(x-z)^{k}Y_W^{*}(v,x)\alpha\nonumber\\
& &-e^{\frac{j_1}{T}\pi i}x_{0}^{-1}\delta\left(\frac{z-x}{-x_{0}}\right)\left(\frac{z-x}{x_{0}}\right)^{j_1\over T}
x^{k+{j_2\over T}}(-z+x)^{k}Y_{\sigma}^{R}(v,x)\alpha\nonumber\\
&=&x_{0}^{-1}\delta\left(\frac{x-z}{x_{0}}\right)
x_{0}^{-{j_1\over T}}\left[x^{k+{j_2\over T}}(x-z)^{k+{j_1\over T}}Y_W^{*}(v,x)\alpha\right]
\nonumber\\
& &-x_{0}^{-1}\delta\left(\frac{z-x}{-x_{0}}\right)
x_{0}^{-{j_1\over T}}\left[e^{(k+{j_1\over T})\pi i}x^{k+{j_2\over T}}(z-x)^{k+{j_1\over T}}Y_{\sigma}^{R}(v,x)\alpha
\right]\nonumber\\
&=&x^{-1}\delta\left(\frac{z+x_{0}}{x}\right)x_{0}^{-{j_1\over T}}\left[x^{k+{j_2\over T}}
(x-z)^{k+{j_1\over T}}Y_W^{*}(v,x)\alpha\right]\nonumber\\
&=&x^{-1}\delta\left(\frac{z+x_{0}}{x}\right)x_{0}^{-{j_1\over T}}
\left[(x_{0}+z)^{k+{j_2\over T}}x_{0}^{k+{j_1\over T}}
Y_W^{*}(v,x_{0}+z)\alpha\right]\nonumber\\
& &\;\;\left(\mbox{with }x \mbox{  replaced by }
x_{0}+z\right)\nonumber\\
&=&x^{-1}\delta\left(\frac{z+x_{0}}{x}\right)
x_{0}^{-{j_1\over T}}\left[
(z+x_{0})^{k+{j_2\over T}}x_{0}^{k+{j_1\over T}}
Y_{\sigma}^{L}(v,x_{0})\alpha\right]\nonumber\\
&=&x^{-1}\delta\left(\frac{z+x_{0}}{x}\right)
\left(\frac{z+x_{0}}{x}\right)^{k+{j_2\over T}}
x^{k+{j_2\over T}}x_{0}^{k}
Y_{\sigma}^{L}(v,x_{0})\alpha.
\end{eqnarray*}
Multiplying by $x^{-k-{j_2\over T}}x_{0}^{-k}$ we get (\ref{sjac}). 
\end{proof}

Next, we combine the twisted $V$-module structures
$Y_{\sigma}^{L}(\cdot,x)$ and $Y_{\sigma}^{R}(\cdot,x)$ on ${\mathfrak{D}}_{\sigma_1,\sigma_2}^{(z)}(W)$
into a twisted $V\otimes V$-module structure. 
With $\sigma_1,\sigma_2$ automorphisms of $V$, 
$\sigma_1\otimes \sigma_2$ is an  automorphism of $V\otimes V$ with $(\sigma_1\otimes \sigma_2)^T=1$.

\bt{ttwisted}
Let $(W,Y_W)$ be a $\sigma$-twisted right $V$-module and let $z$ be a nonzero complex number. 
Then the following commutation relation holds for $u,v\in V$:
\begin{eqnarray}\label{elrrl}
Y_{\sigma}^{L}(u,x_{1})Y_{\sigma}^{R}(v,x_{2})=Y_{\sigma}^{R}(v,x_{2})Y_{\sigma}^{L}(u,x_{1})
\end{eqnarray}
on ${\mathfrak{D}}^{(z)}_{\sigma_1,\sigma_2}(W)$. 
Furthermore, ${\mathfrak{D}}_{\sigma_1,\sigma_2}^{(z)}(W)$ is a $\sigma_1\otimes \sigma_2$-twisted $V\otimes V$-module 
where the vertex operator map $Y^{\rm reg}_{\sigma}(\cdot,x)$ is defined by
\begin{eqnarray}
Y^{\rm reg}_{\sigma}(u\otimes v,x)=Y_{\sigma}^{L}(u,x)Y_{\sigma}^{R}(v,x)\   \   \   \mbox{ for }u,v\in V.
\end{eqnarray}
\et

\begin{proof} It suffices to prove the first assertion. Consider homogeneous vectors 
$u\in V^{(j_{1},j_2)},\  v\in V^{(r_1,r_2)}$ with $j_1,j_2,r_1,r_2\in \Z$. 
Let $\alpha\in {\mathfrak{D}}^{(z)}_{\sigma_1,\sigma_2}(W)$. There exists $k\in {\N}$ such that
\begin{eqnarray}
e^{(k+{r_1\over T})\pi i}(z-x_{1})^{k+{r_1\over T}}Y_{\sigma}^{R}(v,x_{1})\alpha
=(x_{1}-z)^{k+{r_1\over T}}Y_W^{*}(v,x_{1})\alpha.
\end{eqnarray}
In view of Lemma \ref{lbasic1}, we may choose $k$ so that the following also holds:
$$e^{(k+{r_1\over T})\pi i}(z-x_{1})^{k+{r_1\over T}}Y^{R}_{\sigma}(v,x_{1})Y_W^{*}(u,x)\alpha
=(x_{1}-z)^{k+{r_1\over T}}Y_W^{*}(v,x_{1})Y_W^{*}(u,x)\alpha.$$
Furthermore,
since both $Y_{\sigma}^{R}$ and $Y_W^{*}$ satisfy weak commutativity, 
 we may choose $k$ such that the following also hold: 
\begin{eqnarray}
& &(x-x_{1})^{k}Y_W^{*}(u,x)Y_W^{*}(v,x_{1})=(x-x_{1})^{k}Y_W^{*}(v,x_{1})Y_W^{*}(u,x),\\
& &(x-x_{1})^{k}Y_{\sigma}^{R}(u,x)Y_{\sigma}^{R}(v,x_{1})
=(x-x_{1})^{k}Y_{\sigma}^{R}(v,x_{1})Y_{\sigma}^{R}(u,x).
\end{eqnarray}
Then using (\ref{sjac}) (Proposition \ref{prelation}), all the identities above, and Lemma \ref{lbasic1} we obtain
\begin{eqnarray*}\label{elrcom1}
& &x^{-1}\delta\left(\frac{z+x_{0}}{x}\right)
\left(\frac{z+x_{0}}{x}\right)^{{j_2\over T}}
(z-x_{1})^{k+{r_1\over T}}(x-x_{1})^{k}Y_{\sigma}^{L}(u,x_{0})Y_{\sigma}^{R}(v,x_{1})\alpha
\nonumber\\
&=&x_{0}^{-1}\delta\left(\frac{x-z}{x_{0}}\right)
\left(\frac{x-z}{x_{0}}\right)^{j_{1}\over T}
(z-x_{1})^{k+{r_1\over T}}(x-x_{1})^{k}Y_W^{*}(u,x)Y_{\sigma}^{R}(v,x_{1})\alpha
\nonumber\\
& &-e^{{j_{1}\over T}\pi i} x_{0}^{-1}\delta\left(\frac{z-x}{-x_{0}}\right)
\left(\frac{z-x}{x_{0}}\right)^{j_{1}\over T}
(z-x_{1})^{k+{r_1\over T}}(x-x_{1})^{k}Y_{\sigma}^{R}(u,x)Y_{\sigma}^{R}(v,x_{1})\alpha
\nonumber\\
&=&x_{0}^{-1}\delta\left(\frac{x-z}{x_{0}}\right)
\left(\frac{x-z}{x_{0}}\right)^{j_{1}\over T} e^{-(k+{r_1\over T})\pi i}
(x_{1}-z)^{k+{r_1\over T}}(x-x_{1})^{k}Y_W^{*}(u,x)Y_W^{*}(v,x_{1})\alpha
\nonumber\\
& &-e^{{j_{1}\over T}\pi i} x_{0}^{-1}\delta\left(\frac{z-x}{-x_{0}}\right)
\left(\frac{z-x}{x_{0}}\right)^{j_{1}\over T}
(z-x_{1})^{k+{r_1\over T}}(x-x_{1})^{k}Y_{\sigma}^{R}(u,x)Y_{\sigma}^{R}(v,x_{1})\alpha
\nonumber\\
&=&x_{0}^{-1}\delta\left(\frac{x-z}{x_{0}}\right)
\left(\frac{x-z}{x_{0}}\right)^{j_{1}\over T}e^{-(k+{r_1\over T})\pi i}
(x_{1}-z)^{k+{r_1\over T}}(x-x_{1})^{k}Y_W^{*}(v,x_{1})Y_W^{*}(u,x)\alpha
\nonumber\\
& &-e^{{j_{1}\over T}\pi i}x_{0}^{-1}\delta\left(\frac{z-x}{-x_{0}}\right)
\left(\frac{z-x}{x_{0}}\right)^{j_{1}\over T}
(z-x_{1})^{k+{r_1\over T}}(x-x_{1})^{k}Y_{\sigma}^{R}(v,x_{1})Y_{\sigma}^{R}(u,x)\alpha
\nonumber\\
&=&x_{0}^{-1}\delta\left(\frac{x-z}{x_{0}}\right)
\left(\frac{x-z}{x_{0}}\right)^{j_{1}\over T}
(z-x_{1})^{k+{r_1\over T}}(x-x_{1})^{k}Y_{\sigma}^{R}(v,x_{1})Y_W^{*}(u,x)\alpha
\nonumber\\
& &-e^{{j_{1}\over T}\pi i}x_{0}^{-1}\delta\left(\frac{z-x}{-x_{0}}\right)
\left(\frac{z-x}{x_{0}}\right)^{j_{1}\over T}
(z-x_{1})^{k+{r_1\over T}}(x-x_{1})^{k}Y_{\sigma}^{R}(v,x_{1})Y_{\sigma}^{R}(u,x)\alpha
\nonumber\\
&=&x^{-1}\delta\left(\frac{z+x_{0}}{x}\right)\left(\frac{z+x_{0}}{x}\right)^{{j_2\over T}}
(z-x_{1})^{k+{r_1\over T}}(x-x_{1})^{k}Y_{\sigma}^{R}(v,x_{1})Y_{\sigma}^{L}(u,x_{0})\alpha.
\end{eqnarray*}
Then applying $\Res_{x}x^{-{j_{1}\over T}}$ and using (\ref{esub}) we get
\begin{eqnarray}\label{elrcom2}
& &(z+x_{0})^{-{j_{1}\over T}}(z-x_{1})^{k+{r_1\over T}}(z+x_{0}-x_{1})^{k}
Y_{\sigma}^{L}(u,x_{0})Y_{\sigma}^{R}(v,x_{1})\alpha
\nonumber\\
&=&(z+x_{0})^{-{j_{1}\over T}}
(z-x_{1})^{k+{r_1\over T}}(z+x_{0}-x_{1})^{k}Y_{\sigma}^{R}(v,x_{1})Y_{\sigma}^{L}(u,x_{0})\alpha.
\end{eqnarray}
Multiplying both sides  by 
$$(z+x_{0})^{j_{1}\over T}(z-x_{1})^{-k-{r_1\over T}}(z+x_{0}-x_{1})^{-k},$$
we obtain
$$Y_{\sigma}^{L}(u,x_{0})Y_{\sigma}^{R}(v,x_{1})\alpha
=Y_{\sigma}^{R}(v,x_{1})Y_{\sigma}^{L}(u,x_{0})\alpha,$$
as desired.
\end{proof}

\section{Intertwining operators, $P(z)$-intertwining maps, and  $\mathfrak{D}_{\sigma_1,\sigma_2}^{(z)}(W)$}

In this section, we study the interplay of intertwining operators, $P(z)$-intertwining maps, 
and  $\mathfrak{D}_{\sigma_1,\sigma_2}^{(z)}(W)$.
In particular, we show that the socle of ${\mathfrak{D}}_{\sigma_1,\sigma_1^{-1}}^{(z)}(V)$ has a Peter-Weyl type decomposition
and we present a canonical connection of graded trace functions associated to $\sigma_1$-twisted $V$-modules
with $\mathfrak{D}_{\sigma_1,\sigma_1^{-1}}^{(-1)}(V)$.

First of all, for any vector space $U$, following \cite{flm} and \cite{fhl}, set
\begin{eqnarray}
U\{x\}=\left\{A(x)=\sum_{h\in {\C}}a(h)x^{h}\;|\; a(h)\in U
\;\;\;\mbox{ for }h\in {\C}\right\}.
\end{eqnarray}
Define $U\{x\}^o$ to be the subspace of $U\{x\}$, 
consisting of $A(x)=\sum_{\lambda\in \C}a(\lambda)x^{\lambda}$ such that 
for any $\lambda\in \C$, $a(\lambda+n)=0$ for all sufficiently negative integers $n$.

Throughout this section, we assume that $V$ is a vertex operator algebra.
Let $g_1,g_2,g_3$ be mutually commuting finite-order automorphisms of $V$ and 
let $T$ be a positive integer such that  $\sigma_k^T=1$ for $k=1,2,3$. 
For $j_1,j_2\in \Z$, set
\begin{eqnarray}
V^{(j_1,j_2)}=\{v\in V\  |  \  g_kv=e^{2\pi ij_k/T}v,\ k=1,2\}.
\end{eqnarray}
Then 
$$V=\bigoplus_{0\le j_1,j_2<T}V^{(j_1,j_2)}.$$
The following definition is due to Xu (see \cite{xu}; cf. \cite{dl}, \cite{fhl}):

\bd{intertwining-operator}
{\em Let $(M_k,Y_{M_k})$ be a weak $g_k$-twisted $V$-module for $k=1,2,3$.
An {\em intertwining operator of type} $\binom{M_3}{M_1\; M_2}$
 is a linear map 
\begin{eqnarray}
 {\Y}(\cdot,x): &&M_1 \to (\mbox{Hom}(M_2,M_3))\{x\} \nonumber\\
&&w\mapsto {\Y}(w,x)=\displaystyle{\sum_{h\in{\C}}w_hx^{-h-1}}
\end{eqnarray}
such that for any $w^{(i)}\in M_i$ ($i=1,2$) and for any $h\in \C$, 
\begin{equation}\label{a1}
w^{(1)}_{h+n}w^{(2)}=0\   \   \mbox{ for $n\in \Q$ sufficiently large},
\end{equation}
 for $v\in V^{(j_1,j_2)}$ with $j_1,j_2\in \Z$ and for $w^{(1)}\in M_1,$ 
\begin{eqnarray}\label{11.rjac}
&&x^{-1}_0\left(\frac{x_1-x_2}{x_0}\right)^{\frac{j_1}{T}}
\delta\left(\frac{x_1-x_2}{x_0}\right)Y_{M_3}(v,x_1){\Y}(w^{(1)},x_2)\nonumber\\
&&
\  \   \    \ -e^{\frac{j_1}{T}\pi i}x_0^{-1}\left(\frac{x_2-x_1}{x_0}\right)^{\frac{j_1}{T}}
\delta\left(\frac{x_2-x_1}{-x_0}\right)
{\Y}(w^{(1)},x_2)Y_{M_2}(v,x_1)\nonumber\\
&=&x_2^{-1}\left(\frac{x_1-x_0}{x_2}\right)^{-\frac{j_2}{T}}\delta
\left(\frac{x_1-x_0}{x_2}\right){\Y}(Y_{M_1}(v,x_0)w^{(1)},x_2)
\end{eqnarray}
 holds on $M_2$, and
\begin{equation}\label{intertwining-derivative}
{\Y}(L(-1)w^{(1)},x)=\frac{d}{dx}{\Y}(w^{(1)},x)\  \   \   \mbox{ for }w^{(1)}\in M_1.
\end{equation}
Denote by $I{M_{3}\choose M_{1}M_{2}}$ the space of intertwining operators of the indicated type.}
\ed

{\bf A subtlety!}  For $\alpha\in \C$, we have
\begin{eqnarray*}
&&\left(\frac{-x_2+x_1}{x_0}\right)^{\alpha}=x_0^{-\alpha}(-x_2+x_1)^{\alpha}=e^{\alpha\pi i}x_0^{-\alpha}(x_2-x_1)^{\alpha},\\
&&\left(\frac{x_2-x_1}{-x_0}\right)^{\alpha}=(-x_0)^{-\alpha}(x_2-x_1)^{\alpha}=e^{-\alpha\pi i}x_0^{-\alpha}(x_2-x_1)^{\alpha},
\end{eqnarray*}
where  $e^{\alpha\pi i}\ne e^{-\alpha\pi i}$ if $\alpha\notin \Z$.

\br{r-compat}
{\em 
For any intertwining operator $\Y$ of type $\binom{M_3}{M_1\; M_2}$ we have
\begin{eqnarray}
{\Y}(w^{(1)},x_2)w^{(2)}\in M_3\{x_2\}^{o}\quad \text{ for }w^{(1)}\in M_1,\ w^{(2)}\in M_2.
\end{eqnarray}
On the other hand, for $v\in V^{(j_1,j_2)}$ with $j_1,j_2\in \Z$, from (\ref{11.rjac}) we have
\begin{eqnarray}
x^{\frac{j_1+j_2}{T}}Y_{M_3}(v,x){\Y}(w^{(1)},x_2)w^{(2)}\in \left(M_3((x))\right)\!\{x_2\}^{o}
\end{eqnarray}
for $w^{(1)}\in M_1,\ w^{(2)}\in M_2$. If $\Y(\cdot,x)$ is nonzero,  it follows that $g_3=g_1g_2$ (see \cite{xu}).}
\er

For convenience, we formulate the following technical lemma (cf. \cite{li-local}):

\bl{lemm-new-delta}
Let $U$ be any vector space, let $\alpha,\beta\in \mathbb{C}$, and let
\begin{eqnarray*}
&&A(x_1,x_2)\in x_1^{-\alpha-\beta}(U((x_1)))\{x_2\}^o,\  \
B(x_1,x_2)\in x_1^{-\beta}(U\{x_2\}^o)((x_1)),\\
&&\quad\quad C(x_0,x_2)\in x_0^{-\alpha}(U\{x_2\}^o)((x_0)).
\end{eqnarray*}
Then 
\begin{eqnarray}\label{formal-jac}
&&x^{-1}_0
\delta\left(\frac{x_1-x_2}{x_0}\right)\left(\frac{x_1-x_2}{x_0}\right)^{\alpha}A(x_1,x_2)\nonumber\\
&&
\  \   \    \ - e^{\alpha\pi i}x_0^{-1}\delta\left(\frac{x_2-x_1}{-x_0}\right)\left(\frac{x_2-x_1}{x_0}\right)^{\alpha}
B(x_1,x_2)\nonumber\\
&=&x_2^{-1}\delta\left(\frac{x_1-x_0}{x_2}\right)\left(\frac{x_1-x_0}{x_2}\right)^{-\beta} C(x_0,x_2)
\end{eqnarray}
 holds if and only if there exist nonnegative integers $k,l$ such that
 \begin{eqnarray}
&& (x_{1}-x_{2})^{k+\alpha}A(x_1,x_2)=e^{(k+\alpha)\pi i}(x_{2}-x_{1})^{k+\alpha}B(x_1,x_2),\label{generalized-comm-formal}\\
&&(x_{0}+x_{2})^{l+\alpha+\beta}A(x_0+x_2,x_2)=(x_{2}+x_{0})^{l+\alpha+\beta}C(x_0,x_2).\label{generalized-assoc-formal}
 \end{eqnarray}
\el

\begin{proof} By assumption, there exist nonnegative integers $k,l$ such that 
$$x_0^{k+\alpha}C(x_0,x_2)\in (U\{x_2\}^o)[[x_0]]\  \text{ and }\  x_1^{l+\beta}B(x_1,x_2)\in (U\{x_2\}^o)[[x_1]].$$
We see that identity (\ref{formal-jac}) immediately implies (\ref{generalized-comm-formal}) and (\ref{generalized-assoc-formal}).
For the converse, assume that (\ref{generalized-comm-formal}) and (\ref{generalized-assoc-formal}) hold. 
Since $x_1^{l+\beta}B(x_1,x_2)\in (U\{x_2\}^o)[[x_1]]$, using (\ref{generalized-comm-formal}) we get
$$(x_1-x_2)^{\alpha+k}x_1^{l+\beta}A(x_1,x_2)\in (U[[x_1]])\{x_2\}^o,$$
so that
$$\left((x_1-x_2)^{\alpha+k}x_1^{l+\beta}A(x_1,x_2)\right)|_{x_1=x_0+x_2}
=\left((x_1-x_2)^{\alpha+k}x_1^{l+\beta}A(x_1,x_2)\right)|_{x_1=x_2+x_0}.$$
Then we have
\begin{eqnarray*}
&&x^{-1}_0
\delta\left(\frac{x_1-x_2}{x_0}\right)\left(\frac{x_1-x_2}{x_0}\right)^{\alpha}x_0^{k+\alpha}x_1^{l+\beta}A(x_1,x_2)\nonumber\\
&&
\  \   \    \ -e^{\alpha\pi i}x_0^{-1}\delta\left(\frac{x_2-x_1}{-x_0}\right)\left(\frac{x_2-x_1}{x_0}\right)^{\alpha}x_0^{k+\alpha}
x_1^{l+\beta}B(x_1,x_2)\nonumber\\
&=&x^{-1}_0
\delta\left(\frac{x_1-x_2}{x_0}\right)(x_1-x_2)^{\alpha+k}x_1^{l+\beta}A(x_1,x_2)\nonumber\\
&&
\  \   \    \ -x_0^{-1}\delta\left(\frac{x_2-x_1}{-x_0}\right)e^{(\alpha+k)\pi i}(x_2-x_1)^{\alpha+k}x_1^{l+\beta}B(x_1,x_2)\nonumber\\
&=&x_1^{-1}\delta\left(\frac{x_2+x_0}{x_1}\right)\left((x_1-x_2)^{\alpha+k}x_1^{l+\beta}A(x_1,x_2)\right)\nonumber\\
&=&x_1^{-1}\delta\left(\frac{x_2+x_0}{x_1}\right)\left((x_1-x_2)^{\alpha+k}x_1^{l+\beta}A(x_1,x_2)\right)|_{x_1=x_2+x_0}\nonumber\\
&=&x_1^{-1}\delta\left(\frac{x_2+x_0}{x_1}\right)\left((x_1-x_2)^{\alpha+k}x_1^{l+\beta}A(x_1,x_2)\right)|_{x_1=x_0+x_2}\nonumber\\
&=&x_1^{-1}\delta\left(\frac{x_2+x_0}{x_1}\right)\left(x_0^{\alpha+k}(x_0+x_2)^{l+\beta}A(x_0+x_2,x_2)\right)\nonumber\\
&=&x_1^{-1}\delta\left(\frac{x_2+x_0}{x_1}\right)\left(x_0^{\alpha+k}(x_2+x_0)^{l+\beta}C(x_0,x_2)\right)\nonumber\\
&=&x_1^{-1}\delta\left(\frac{x_2+x_0}{x_1}\right)\left(\frac{x_2+x_0}{x_1}\right)^{l+\beta}
\left(x_0^{\alpha+k}x_1^{l+\beta}C(x_0,x_2)\right)\nonumber\\
&=&x_1^{-1}\delta\left(\frac{x_2+x_0}{x_1}\right)\left(\frac{x_2+x_0}{x_1}\right)^{\beta} x_0^{k+\alpha}x_1^{l+\beta}C(x_0,x_2),
\end{eqnarray*}
which readily implies (\ref{formal-jac}).
\end{proof}

With Lemma \ref{lemm-new-delta} we immediately have:

\bl{rtwistedcommassoc} 
The generalized Jacobi identity (\ref{11.rjac}) is equivalent to the 
following properties:
(1) For $v\in V^{(j_1,j_2)}\  (j_1,j_2\in \Z),\ w\in M_1$, there exists $k\in {\N}$ such that
\begin{eqnarray}\label{etwistedweakcomm}
(x_{1}-x_{2})^{k+{j_1\over T}}Y_{M_3}(v,x_{1}){\Y}(w,x_{2})
=e^{(k+{j_1\over T})\pi i}(x_{2}-x_{1})^{k+{j_1\over T}}{\Y}(w,x_{2})Y_{M_2}(v,x_{1}).
\end{eqnarray}
(2) For $v\in V^{(j_1,j_2)}\  (j_1,j_2\in \Z),\ w\in M_1,\  w^2\in M_2$, there exists $l\in {\N}$ such that
\begin{eqnarray}\label{etwistedweakassoc}
(x_{0}+x_{2})^{l+{j_1+j_2\over T}}Y_{M_3}(v,x_{0}+x_{2}){\Y}(w,x_{2})w^2
=(x_{2}+x_{0})^{l+{j_1+j_2\over T}}{\Y}(Y_{M_1}(v,x_0)w,x_2)w^2.
\end{eqnarray}
\el

Note that for $\mu\in \mathbb{C}$, we have (see \cite{dl})
\begin{eqnarray*}
z^{-1}\delta\left(\frac{x_{1}-x_{0}}{z}\right)\left(\frac{x_{1}-x_{0}}{z}\right)^{\mu}
=x_1^{-1}\delta\left(\frac{z+x_{0}}{x_1}\right)\left(\frac{z+x_{0}}{x_1}\right)^{-\mu},
\end{eqnarray*}
$$x_{0}^{-1}\delta\left(\frac{x_{1}-z}{x_{0}}\right)\left(\frac{x_{1}-z}{x_{0}}\right)^{\mu}
=x_{1}^{-1}\delta\left(\frac{x_{0}+z}{x_{1}}\right)\left(\frac{x_{0}+z}{x_{1}}\right)^{-\mu}.$$

Let $\sigma_1,\sigma_2, \sigma$, and $T$ be  given as before, in particular,
$\sigma=\sigma_1^{-1}\sigma_2^{-1}$.
Let $z$ be a nonzero complex number.
The following is a twisted analogue of the notion of $P(z)$-intertwining map due to \cite{hl1}: 
 
\bd{def-p(z)im} 
{\em Let $W$ be a weak $\sigma$-twisted $V$-module, $W_1$ a weak $\sigma_1$-twisted $V$-module, and
$W_2$ a weak $\sigma_2$-twisted $V$-module.
A {\em $P(z)$-intertwining map} of type ${W^{*}\choose W_{1}W_{2}}$ is a linear map
$F: W_{1}\otimes W_{2}\rightarrow W^{*}$ such that 
for $v\in V^{(j_1,j_2)}, \ j_1,j_2\in \Z,\    w_{(1)}\in W_{1},\   w_{(2)}\in W_{(2)}$,
\begin{eqnarray}\label{P(z)-Jacobi}
& &x_{0}^{-1}\delta\left(\frac{x_{1}-z}{x_{0}}\right)\left(\frac{x_{1}-z}{x_{0}}\right)^{\frac{j_1}{T}}
Y_W^{*}(v,x_{1})F(w_{(1)}\otimes w_{(2)})\nonumber\\
&&\hspace{0.5cm}  - e^{\frac{j_1}{T}\pi i}x_{0}^{-1}\delta\left(\frac{z-x_{1}}{-x_{0}}\right)\left(\frac{z-x_{1}}{x_{0}}\right)^{\frac{j_1}{T}}
F(w_{(1)}\otimes Y(v,x_{1})w_{(2)})\nonumber\\
&= &z^{-1}\delta\left(\frac{x_{1}-x_{0}}{z}\right)\left(\frac{x_{1}-x_{0}}{z}\right)^{-\frac{j_2}{T}}
F(Y(v,x_{0})w_{(1)}\otimes w_{(2)}).
\end{eqnarray}}
\ed

By a routine argument we get:

\bl{P(z)-map-twistedcommassoc} 
A linear map $F: W_1\otimes W_2\rightarrow W^{*}$ is a $P(z)$-intertwining map if and only if the following two conditions hold
for $v\in V^{(j_1,j_2)}\  (j_1,j_2\in \Z),\ w^{(1)}\in W_1,\ w^{(2)}\in W_2$:\\
 (1) There exists $k\in {\N}$ such that
\begin{eqnarray}\label{etwistedweakcomm}
(x-z)^{k+{j_1\over T}}Y_{W}^{*}(v,x)F(w^{(1)}\otimes w^{(2)})
=e^{(k+{j_1\over T})\pi i}(z-x)^{k+{j_1\over T}}F(w^{(1)}\otimes Y(v,x)w^{(2)}).
\end{eqnarray}
(2) There exists $l\in {\N}$ such that
\begin{eqnarray}\label{etwistedweakassoc}
(x+z)^{l+{j_1+j_2\over T}}Y_{W}^{*}(v,x+z)F(w^{(1)}\otimes w^{(2)})
=(z+x)^{l+{j_1+j_2\over T}}F(Y(v,x)w^{(1)}\otimes w^{(2)}).
\end{eqnarray}
\el

As the first main result of this section, we have:

\bt{t-Int-map-dual} 
Let $W$ be a weak $\sigma$-twisted $V$-module, $W_1$ a weak $\sigma_1$-twisted $V$-module, and
$W_2$ a weak $\sigma_2$-twisted $V$-module.
 A linear map $F: W_1\otimes W_2\rightarrow W^{*}$ is a $P(z)$-intertwining map if and only if 
it is a  homomorphism of weak $\sigma_1\times \sigma_2$-twisted $V\otimes V$-modules into 
$\mathfrak{D}_{\sigma_1,\sigma_2}^{(z)}(W)$.
\et

\begin{proof} Assume that $F: W_1\otimes W_2\rightarrow W^{*}$ is a $P(z)$-intertwining map. We must prove 
$$F(w_{(1)}\otimes w_{(2)})\in \mathfrak{D}_{\sigma_1,\sigma_2}^{(z)}(W)
\quad \text{for all }w_{(1)}\in W_{1},\ w_{(2)}\in W_{2}, $$
and for every $v\in V$,
\begin{eqnarray}
& &Y_{\sigma}^{R}(v,x)F(w_{(1)}\otimes w_{(2)})
=F(w_{(1)}\otimes Y(v,x)w_{(2)}),\label{eYRF}\\
& &Y_{\sigma}^{L}(v,x)F(w_{(1)}\otimes w_{(2)})
=F(Y(v,x)w_{(1)}\otimes w_{(2)}).\label{eYLF}
\end{eqnarray}
Assume $v\in V^{(j_1,j_2)}\ (j_1,j_2\in \Z), \  w_{(1)}\in W_{1},\  w_{(2)}\in W_{2}$. 
As $x_0^{\frac{j_1}{T}}Y(v,x_0)w_{(1)}\in W_1((x_0))$, there exists $k\in \N$ 
such that $x_0^{k+\frac{j_1}{T}}Y(v,x_0)w_{(1)}\in W_1[[x_0]]$.
From (\ref{P(z)-Jacobi}), we get 
\begin{eqnarray}\label{e-glocality-proof}
(x-z)^{k+{j_1\over T}}Y_W^{*}(v,x)
F(w_{(1)}\otimes w_{(2)})=e^{(k+{j_1\over T})\pi i}(z-x)^{k+{j_1\over T}}F(w_{(1)}\otimes Y(v,x)w_{(2)}).
\end{eqnarray}
Noticing that $x^{\frac{j_2}{T}}Y(v,x)w_{(2)}\in W_2((x))$,  we conclude
$$x^{\frac{j_2}{T}}(x-z)^{k+{j_1\over T}}Y_W^{*}(v,x)
F(w_{(1)}\otimes w_{(2)})\in W^{*}((x)).$$
Thus $F(w_{(1)}\otimes w_{(2)})\in \mathfrak{D}_{\sigma_1,\sigma_2}^{(z)}(W)$.
Furthermore, by Corollary \ref{cbasic} and (\ref{e-glocality-proof}), we have
\begin{eqnarray}\label{etemp}
& &e^{(k+{j_1\over T})\pi i}(z-x)^{k+{j_1\over T}}Y_{\sigma}^{R}(v,x)
F(w_{(1)}\otimes w_{(2)})\nonumber\\
&=&(x-z)^{k+{j_1\over T}}Y_W^{*}(v,x)
F(w_{(1)}\otimes w_{(2)})
\nonumber\\
&=&e^{(k+{j_1\over T})\pi i}(z-x)^{k+{j_1\over T}}F(w_{(1)}\otimes Y(v,x)w_{(2)}).
\end{eqnarray}
As
$$Y_{\sigma}^{R}(v,x)F(w_{(1)}\otimes w_{(2)}),\;\;\;\; 
F(w_{(1)}\otimes Y(v,x)w_{(2)})
\in x^{-{j_2\over T}}{\mathfrak{D}}_{\sigma_1,\sigma_2}^{(z)}(W)((x)),$$
multiplying both sides of (\ref{etemp}) by $e^{-(k+{j_1\over T})\pi i}(z-x)^{-k-{j_1\over T}}$  we obtain (\ref{eYRF}).

Similarly, from (\ref{P(z)-Jacobi}) there exists $l\in {\N}$ such that
$$(x+z)^{l+{j_2\over T}}Y_W^{*}(v, x+z)F(w_{(1)}\otimes w_{(2)})=
(z+x)^{l+{j_2\over T}}F(Y(v,x)w_{(1)}\otimes w_{(2)}),$$
which  lies in $x^{-{j_1\over T}}W^{*}((x))$.
Then by Definition \ref{dleftaction} we have
\begin{eqnarray}\label{etemp-l}
& &(z+x)^{l+{j_2\over T}}Y_{\sigma}^{L}(v,x)
F(w_{(1)}\otimes w_{(2)})\nonumber\\
&=&(x+z)^{l+{j_2\over T}}
Y_W^{*}(v,x+z)F(w_{(1)}\otimes w_{(2)})\nonumber\\
&=&(z+x)^{l+{j_2\over T}}F(Y(v,x)w_{(1)}\otimes w_{(2)}).
\end{eqnarray}
Again, as both
$Y_{\sigma}^{L}(v,x)F(w_{(1)}\otimes w_{(2)})$ and  $F(Y(v,x)w_{(1)}\otimes w_{(2)})$
lie in $x^{-{j_1\over T}}{\mathfrak{D}}_{\sigma_1,\sigma_2}^{(z)}(W)((x))$,
we can multiply (\ref{etemp-l}) by $(z+x)^{-l-{j_2\over T}}$ to obtain (\ref{eYLF}). 

For the converse, assume that $F$ is a homomorphism of $\sigma_1\times \sigma_2$-twisted $V\otimes V$-modules
 into $\mathfrak{D}_{\sigma_1,\sigma_2}^{(z)}(W)$. Let $v\in V^{(j_1,j_2)}\ (j_1,j_2\in \Z), \  w_{(1)}\in W_{1},\  w_{(2)}\in W_{2}$. 
 As $F(w_{(1)}\otimes w_{(2)})\in \mathfrak{D}_{\sigma_1,\sigma_2}^{(z)}(W)$, there exists $k\in \N$ such that 
 $$x^{\frac{j_2}{T}}(x-z)^{k+{j_1\over T}}Y_W^{*}(v,x)
F(w_{(1)}\otimes w_{(2)})\in W^{*}((x)).$$
Then 
\begin{eqnarray}\label{etemp-converse}
&&(x-z)^{k+{j_1\over T}}Y_W^{*}(v,x)
F(w_{(1)}\otimes w_{(2)})
\nonumber\\
&=&e^{(k+{j_1\over T})\pi i}(z-x)^{k+{j_1\over T}}Y_{\sigma}^{R}(v,x)
F(w_{(1)}\otimes w_{(2)})\nonumber\\
&=&e^{(k+{j_1\over T})\pi i}(z-x)^{k+{j_1\over T}}F(w_{(1)}\otimes Y(v,x)w_{(2)}).
\end{eqnarray}
On the other hand, by Lemma \ref{pleftrightcondition} there exists $l\in \N$ such that 
$$(x+z)^{l+{j_2\over T}}Y_W^{*}(v,x+z)
F(w_{(1)}\otimes w_{(2)})\in x^{-\frac{j_1}{T}}W^{*}((x)).$$
Then by Definition \ref{dleftaction} we get
\begin{eqnarray}\label{etemp-converse-assoc}
(x+z)^{l+{j_2\over T}}Y_W^{*}(v,x+z)
F(w_{(1)}\otimes w_{(2)})&=&(z+x)^{l+{j_2\over T}}Y_{\sigma}^{L}(v,x)F(w_{(1)}\otimes w_{(2)})\nonumber\\
&=&(z+x)^{l+{j_2\over T}}F(Y(v,x)w_{(1)}\otimes w_{(2)}).
\end{eqnarray}
With (\ref{etemp-converse}) and (\ref{etemp-converse-assoc}), by Lemma \ref{P(z)-map-twistedcommassoc} 
 we conclude that $F$ is a $P(z)$-intertwining map.
Now, the proof is complete.
\end{proof}

Denote by $I_{P(z)}{W^{*}\choose W_{1}\; W_{2}}$ the space of $P(z)$-intertwining maps of the indicated type.
The following follows from Theorem \ref{t-Int-map-dual} immediately:

\bc{socle}
Let ${\rm soc}({\mathfrak{D}}_{\sigma_1,\sigma_2}^{(z)}(W))$ denote the sum of all irreducible $\sigma_1\times \sigma_2$-twisted
$V\otimes V$-submodules of ${\mathfrak{D}}_{\sigma_1,\sigma_2}^{(z)}(W)$, called the socle. Then
\begin{eqnarray}
{\rm soc}({\mathfrak{D}}_{\sigma_1,\sigma_2}^{(z)}(W))= \coprod_{W_{1}\in {\rm Irr}_{\sigma_1}(V),W_2\in {\rm Irr}_{\sigma_2}(V)}
I_{P(z)}{W^{*}\choose W_{1}\; W_{2}}\otimes (W_{1}\otimes W_{2})
\end{eqnarray}
as (weak) $\sigma_1\times \sigma_2$-twisted $V\otimes V$-modules. 
\ec

In the following, we slightly extend a result of \cite{hl2}, giving an equivalence between the notions of 
intertwining operator and $P(z)$-intertwining map.
We assume that $W$ is a weak $\sigma$-twisted $V$-module on which $L(0)$ is semisimple. 
Then $$W=\oplus_{h\in \C}W_{(h)}, \  \  \text{where }W_{(h)}=\{ w\in W\ |\ L(0)w=hw\}.$$
Set 
$$W'=\oplus_{h\in \C}W_{(h)}^{*}\subset W^{*}.$$

Define a linear map
$$Y_W^{*}(\cdot,x):\ V\rightarrow (\End W^{*})[[x^{\frac{1}{T}},x^{-\frac{1}{T}}]]$$
as before by
\begin{eqnarray}\label{Y*-Yo}
\< Y_W^{*}(v,x)\alpha,w\>=\< \alpha,Y_W^{o}(v,x)w\>
\end{eqnarray}
for $v\in V,\ \alpha\in W^{*},\ w\in W$.

For $v\in V$, set
\begin{eqnarray}
Y_W^{*}(v,x)=\sum_{n\in \Z}v_{n}^{*}x^{-n-1}\ \  (\text{where }v_n^{*}\in \End W^{*}).
\end{eqnarray}

\br{rem-full-dual-virasoro}
{\em If $W$ satisfies the condition that for every $h\in \C$, $W_{(h+\frac{n}{T})}=0$
for all sufficiently negative integers $n$, then $W'$ is a weak $\sigma^{-1}$-twisted $V$-module.
In general,  $W'$ is not a weak $\sigma^{-1}$-twisted $V$-module.
On the other hand, $W^{*}$ is always a module for the Virasoro algebra with $L(n)$ acting as $\omega^{*}_{n+1}$ 
for $n\in \Z$ of the same central charge $c_V$ as for $V$ and $W$, where
\begin{eqnarray}
\langle L(n)\alpha, w\rangle=\langle \alpha, L(-n)w\rangle \quad \text{for }n\in \Z,\ \alpha\in W^{*},\ w\in W. 
\end{eqnarray}}
\er

\bd{def-W-star}
{\em Denote by $W^{\star}$ the subspace of $W^{*}$, consisting of all $f\in W^{*}$ such that
 for every $h\in \C$,  $\<f, W_{(h+\frac{n}{T})}\>=0$ for all sufficiently negative integers $n$.}
 \ed

\bl{lem-in-Wstar}
Let $f\in W^{\star}$. Then for any $v\in V,\ q\in \frac{1}{T}\Z,\ \lambda_r\in \C$ with $r\in \N$,
\begin{eqnarray}
\sum_{r\in \N}\lambda_r v^{*}_{q-r}f\in W^{\star}.
\end{eqnarray}
Furthermore,  $W^{\star}\cap {\mathfrak{D}}_{\sigma_1,\sigma_2}^{(z)}(W)$ 
is a weak $\sigma_1\times \sigma_2$-twisted $V\otimes V$-submodule of 
${\mathfrak{D}}_{\sigma_1,\sigma_2}^{(z)}(W)$.
\el

\begin{proof} Let $u\in V$ be homogeneous and let $q\in \frac{1}{T}\Z,\ h\in \C$. Since
$$u_{q+r}W_{(h+\frac{n}{T})}\subset W_{(h+\frac{n}{T}+\wt u-q-r-1)} $$
for $r, n\in \Z$ and since $f\in W^{\star}$, there exists $k\in \Z$ (depending on $u,q,h$) such that 
$$\< f, u_{q+r}W_{(h+\frac{n}{T})}\>=0\quad  \text{for all integers }r\ge 0,\ n\le k.$$ 
Now, let $v\in V^{(j_1,j_2)}$ be any homogeneous vector with $j_1,j_2\in \Z$. Write
 $$e^{xL(1)}(-x^{-2})^{L(0)}v=x^{p_1}u^{(1)}+\cdots +x^{p_s}u^{(s)},$$
where $p_1,\dots,p_s\in \Z$, and $u^{(1)},\dots,u^{(s)}\in V$ are homogeneous.  
Then for any fixed $q\in \frac{1}{T}\Z,\ h\in \C$, there exists $k'\in \Z$ such that 
\begin{eqnarray}
\<v_{q-r}^{*}f,W_{(h+\frac{n}{T})}\>=\sum_{i=1}^s\<f,u^{(i)}_{-p_i-q-2+r}W_{(h+\frac{n}{T})}\>=0
\end{eqnarray}
for all $r\ge 0,\ n\le k'$. This in particular proves that $v^{*}_qf\in W^{\star}$. 
Note that for every $w\in W$, $\<v_{q-r}^{*}f,w\>=\sum_{i=1}^s\<f,u^{(i)}_{-p_i-q-2+r}w\rangle=0$ for $r$ sufficiently positive. Thus
any infinite sum $\sum_{r\in \N}\lambda_r v^{*}_{q-r}f$ with $\lambda_r\in \C$ exists in $W^{*}$ and lies in $W^{\star}$.
This proves the first assertion. 

Now, assume $v\in V^{(j_1,j_2)},\ f\in W^{\star}\cap {\mathfrak{D}}_{\sigma_1,\sigma_2}^{(z)}(W)$.
By the first assertion, for any $\beta\in \C,\ p, q\in \Z$, 
$$\Res_{x}x^{q-\beta}(x-z)^{\beta+p}Y_W^{*}(v,x)f=\sum_{r\ge 0}\binom{\beta+p}{r}(-z)^rv^{*}_{p+q-r}f\in W^{\star},$$ 
which implies
\begin{eqnarray}\label{beta-W-star}
x^{-\beta}(x-z)^{\beta+p}Y_W^{*}(v,x)f\in W^{\star}[[x^{\frac{1}{T}},x^{-\frac{1}{T}}]].
\end{eqnarray}
Then from the definition of $Y_{\sigma}^{R}(v,x)f$ we get 
$x^{\frac{j_2}{T}}Y_{\sigma}^{R}(v,x)f\in (W^{\star}\cap {\mathfrak{D}}_{\sigma_1,\sigma_2}^{(z)}(W))((x))$. 
 Furthermore, using Proposition \ref{prelation} we obtain 
$x^{\frac{j_1}{T}}Y_{\sigma}^{L}(v,x)f\in (W^{\star}\cap {\mathfrak{D}}_{\sigma_1,\sigma_2}^{(z)}(W))((x))$. 
This proves the second assertion, completing the proof.
\end{proof}

\br{rem-Wstar}
{\em By Lemma \ref{lem-in-Wstar}, $W^{\star}$ is closed under the actions of $v_n^{*}$ for $v\in V,\ n\in \frac{1}{T}\Z$.
 This in particular implies that $W^{\star}$ is a submodule of $W^{*}$ for the Virasoro algebra.}
\er

In addition to $W$, we also assume that 
$W_1$ is a weak $\sigma_1$-twisted $V$-module and $W_2$
is a weak $\sigma_2$-twisted $V$-module such that $L(0)$ is semisimple on $W_1$ and $W_2$.

\br{rintL(-1)}
{\em Even though $W'$ is not a weak $\sigma^{-1}$-twisted $V$-module in general, Definition \ref{intertwining-operator} 
for an intertwining operator of type $\binom{W'}{W_1\; W_2}$ works fine.
Then we define an intertwining operator of type $\binom{W'}{W_1\; W_2}$ this way.}
\er

Suppose $\mathcal{Y}(\cdot,x)$ is an intertwining operator of type $\binom{W'}{W_1\; W_2}$. 
Let $w_{(1)}\in W_{1}$ and $w_{(2)}\in W_{2}$ be homogeneous. 
From Jacobi identity (\ref{P(z)-Jacobi}) we have (see \cite{fhl})
\begin{eqnarray}\label{L(0)-bracket}
[L(0),\Y(w_{(1)},x)]w_{(2)}=x\Y(L(-1)w_{(1)},x)w_{(2)}+\Y(L(0)w_{(1)},x)w_{(2)}.
\end{eqnarray}
This particularly implies 
\begin{eqnarray}\label{w1nw2}
L(0)(w_{(1)})_{n}w_{(2)}=(\wt w_{(1)}+\wt w_{(2)}-n-1)(w_{(1)})_{n}w_{(2)}
\end{eqnarray}
for $n\in \C$. From \cite{hl1},
 (\ref{L(0)-bracket}) and (\ref{w1nw2}) also imply the $L(-1)$-derivative property (\ref{intertwining-derivative}).
 
Recall that 
$$\log z=\ln |z|+i{\rm arg}(z)\ \text{ with }0\le {\rm arg}(z)<2\pi.$$
Then (see \cite{hl1})
$${\Y}(w_{(1)},e^{\log z})w_{(2)}:
={\Y}(w_{(1)},x)w_{(2)}|_{x=e^{\log z}}\in W^{*}$$
and
\begin{eqnarray}
(Y_{W'}(v,x){\Y}(w_{(1)},y)w_{(2)})|_{y=e^{\log z}}
=Y_W^{*}(v,x)({\Y}(w_{(1)},e^{\log z})w_{(2)})
\end{eqnarray}
for $v\in V,\  w_{(1)}\in W_{1},\  w_{(2)}\in W_{2}$.
Let $I_{P(z)}{W^{\star}\choose W_{1}\; W_{2}}$ denote the space of $P(z)$-intertwining maps of the indicated type.
With  (\ref{w1nw2}) and (\ref{a1}), by Theorem \ref{t-Int-map-dual} we immediately have:

\bp{psub}
Let $\Y(\cdot,x)$ be an intertwining operator of type
${W'\choose W_{1}\; W_{2}}$. Then the linear map 
$F_{{\Y}}^{(z)}: W_1\otimes W_2\rightarrow  W^{*}$, defined by
\begin{eqnarray}
F_{{\Y}}^{(z)}(w_{(1)}\otimes w_{(2)})={\Y}(w_{(1)},e^{\log z})w_{(2)}
\end{eqnarray}
for $w_{(1)}\in W_{1},\  w_{(2)}\in W_{2}$, is a $P(z)$-intertwining map of type 
$\binom{W^{*}}{W_1\ W_2}$ such that
\begin{eqnarray}
F_{{\Y}}^{(z)}(W_{1}\otimes W_{2}) \subset W^{\star}\cap {\mathfrak{D}}_{\sigma_1,\sigma_2}^{(z)}(W).
\end{eqnarray}
Furthermore, $F_{{\Y}}^{(z)}$  is a homomorphism of weak $\sigma_1\otimes \sigma_2$-twisted $V\otimes V$-modules 
from $W_{1}\otimes W_{2}$ to $W^{\star}\cap {\mathfrak{D}}_{\sigma_1,\sigma_2}^{(z)}(W)$.
\ep

By Proposition \ref{psub},  we have a linear map
\begin{eqnarray}
F^{(z)}:& & I{W'\choose W_{1}\; W_{2}}\rightarrow 
{\rm Hom}_{V\otimes V}\left(W_{1}\otimes W_{2}, W^{\star}\cap {\mathfrak{D}}_{\sigma_1,\sigma_2}^{(z)}(W)\right)
\nonumber\\
& &{\Y}\mapsto F_{{\Y}}^{(z)}.
\end{eqnarray}
In case $\sigma_1=1=\sigma_2$, it was proved in \cite{hl3} (Proposition 12.2) that
$F^{(z)}$ is a linear isomorphism from the space of intertwining operators to the space of $P(z)$-intertwining maps.
Next, we follow \cite{hl3} to show that this is also true for this generalization.

Recall the following result from \cite{fhl} (Lemma 5.2.3):

\bl{lfhl523}
Let $\Y(\cdot,x)$ be an intertwining operator of type ${W'\choose W_{1}\; W_{2}}$. Then 
\begin{eqnarray}\label{L(0)con}
x^{L(0)}\Y(w_{(1)},x_{0})x^{-L(0)}=\Y(x^{L(0)}w_{(1)},xx_{0})
\end{eqnarray}
for $w_{(1)}\in W_{1}$, where $x$ and $x_{0}$ are independent commuting
formal variables.
\el

Note  that $x^{\pm L(0)}$ are linear maps from $W$ to $W\{x\}$, 
where $x^{\pm L(0)}w=x^{\pm h}w$ for $w\in W_{(h)},\  h\in {\C}$. 
Then view $x^{\pm L(0)}$ as linear maps from $W^{*}$ to $W^{*}\{x\}$ by
\begin{eqnarray}
\<x^{\pm L(0)}f,w\>=\<f, x^{\pm L(0)}w\> \quad \text{ for }f\in W^{*},\ w\in W.  
\end{eqnarray}
The following conjugation formula holds for $v\in V,\   \alpha\in W^{*}$:
\begin{eqnarray}\label{econ*}
x_{1}^{L(0)}Y_W^{*}(v,x)x_{1}^{-L(0)}\alpha=Y_W^{*}(x_{1}^{L(0)}v,x_{1}x)\alpha.
\end{eqnarray}
On the other hand, view $e^{\pm (\log z)L(0)}$ as  linear operators on $W^{*}$ by
\begin{eqnarray}
\<e^{\pm (\log z)L(0)}f,w\>=\<f,e^{\pm (\log z)L(0)}w\>
\end{eqnarray}
for $f\in W^{*},\  w\in W$. Specializing $x_{0}=e^{\log z}$ in  (\ref{L(0)con}) we get
\begin{eqnarray}\label{L(0)conz}
x^{L(0)}\Y(w_{(1)},e^{\log z})x^{-L(0)}w_{(2)}=\Y(x^{L(0)}w_{(1)},xe^{\log z})w_{(2)}
\end{eqnarray}
for $w_{(1)}\in W_1,\  w_{(2)}\in W_{2}$. Then
\begin{eqnarray}\label{recov}
\Y(w_{(1)},x)w_{(2)}=y^{L(0)} \Y\left(y^{-L(0)}w_{(1)},e^{\log z}\right)
y^{-L(0)}w_{(2)}|_{y=xe^{-\log z}}
\end{eqnarray}
(see \cite{hl1}, (4.15)).

We have:

\bt{thm-main} 
Let $W, W_1, W_2$, and $z$ be given as above.
Then the linear map $F^{(z)}$ is a linear isomorphism from $I{W'\choose W_{1}\; W_{2}}$ 
to $I_{P(z)}{W^{\star}\choose W_{1}\; W_{2}}$, which coincides with
${\rm Hom}_{V\otimes V}(W_{1}\otimes W_{2},W^{\star}\cap {\mathfrak{D}}_{\sigma_1,\sigma_2}^{(z)}(W))$.
\et

\begin{proof} The same arguments of (\cite{hl3}, Proposition 12.2) show that 
$F^{(z)}$ is injective. Next, we follow \cite{hl3} to prove that $F^{(z)}$ is also onto.

Let $\psi: W_{1}\otimes W_{2}\rightarrow W^{\star}\cap {\mathfrak{D}}_{\sigma_1,\sigma_2}^{(z)}(W)$
 be a  homomorphism of weak $\sigma_1\otimes \sigma_2$-twisted $V\otimes V$-modules.
 Define a bilinear map $\Y(\cdot,x)\cdot: \ W_1\times W_2\rightarrow W'\{x\}$ by (\ref{recov}).
 For $h\in \C$, let $p_h: W^{*}\rightarrow W_{(h)}^{*}$ be the projection map.
 For $w_{(1)}\in W_{1},\ w_{(2)}\in W_{2},\ h\in \C$, set
 $$\psi(w_{(1)}\otimes w_{(2)})_h=p_{h}(\psi(w_{(1)}\otimes w_{(2)})).$$
Assume $w_{(1)}\in W_1,\ w_{(2)}\in W_2$ are homogeneous with
$\wt w_{(1)}=h_{1},\   \wt w_{(2)}=h_{2}$. Then
\begin{eqnarray}
\Y(w_{(1)},x)w_{(2)}=\sum_{h\in {\C}}x^{h-h_{1}-h_{2}}e^{(h_{1}+h_{2}-h)\log z}
\psi(w_{(1)}\otimes w_{(2)})_{h}.
\end{eqnarray}
Writing
$$\Y(w_{(1)},x)w_{(2)}=\sum_{n\in {\C}}(w_{(1)})_nw_{(2)}x^{-n-1},  $$
we have $(w_{(1)})_nw_{(2)}=e^{(n+1)\log z}\psi(w_{(1)}\otimes w_{(2)})_{h_1+h_2-n-1}$, so that
\begin{eqnarray}\label{weight}
L(0)(w_{(1)})_{n}w_{(2)}=(h_{1}+h_{2}-n-1)(w_{(1)})_{n}w_{(2)}.
\end{eqnarray}
Since $\psi(W_{1}\otimes W_{2})\subset W^{\star}$, we conclude $\Y(w_{(1)},x)w_{(2)}\in W'\{x\}^{o}$.
In the following, we prove that the generalized weak commutativity and associativity relations hold.

Let $v\in V^{(j_1,j_2)},\  w_{(1)}\in W_{1},\  w_{(2)}\in W_{2}$ be homogeneous. 
Since $\psi(w_{(1)}\otimes w_{(2)})\in {\mathfrak{D}}_{\sigma_1,\sigma_2}^{(z)}(W)$, by 
Corollary \ref{cbasic} there exists $k\in {\N}$ such that
$$(x-z)^{k+{j_1\over T}}Y_W^{*}(v,x)\psi(w_{(1)}\otimes w_{(2)})
=e^{(k+{j_1\over T})\pi i}(z-x)^{k+{j_1\over T}}
Y_{\sigma}^{R}(v,x)\psi(w_{(1)}\otimes w_{(2)}),$$
hence
\begin{eqnarray}\label{e4.38}
& &(x_{1}y^{-1}-z)^{k+{j_1\over T}}
Y_W^{*}(y^{-L(0)}v,x_{1}y^{-1})
\psi(y^{-L(0)}w_{(1)}\otimes y^{-L(0)}w_{(2)})\nonumber\\
&=&e^{(k+{j_1\over T})\pi i}(z-x_{1}y^{-1})^{k+{j_1\over T}}
Y_{\sigma}^{R}(y^{-L(0)}v,x_{1}y^{-1})
\psi(y^{-L(0)}w_{(1)}\otimes y^{-L(0)}w_{(2)}).
\end{eqnarray}
Noticing that 
$$(x_{1}y^{-1}-z)^{k+{j_1\over T}}|_{y=x_{2}e^{-\log z}}
=z^{k}e^{{j_1\over T}\log z}x_{2}^{-k-{j_1\over T}}(x_{1}-x_{2})^{k+{j_1\over T}}$$
and
$$(z-x_{1}y^{-1})^{k+{j_1\over T}}|_{y=x_{2}e^{-\log z}}
=z^{k}e^{{j_1\over T}\log z}x_{2}^{-k-{j_1\over T}}(x_{2}-x_{1})^{k+{j_1\over T}},$$
then from  (\ref{e4.38}) we get
\begin{eqnarray*}\label{e437}
& &(x_{1}-x_{2})^{k+{j_1\over T}}
Y_W^{*}(y^{-L(0)}v,x_{1}y^{-1})
\psi(y^{-L(0)}w_{(1)}\otimes y^{-L(0)}w_{(2)})|_{y=x_{2}e^{-\log z}}
\nonumber\\
&=&e^{(k+{j_1\over T})\pi i}(x_{2}-x_{1})^{k+{j_1\over T}}
Y_{\sigma}^{R}(y^{-L(0)}v,x_{1}y^{-1})\psi(y^{-L(0)}w_{(1)}\otimes 
y^{-L(0)}w_{(2)})|_{y=x_{2}e^{-\log z}}.
\end{eqnarray*}
Using this and (\ref{econ*}) we obtain the desired generalized weak commutativity relation as
\begin{eqnarray*}
& &(x_{1}-x_{2})^{k+{j_1\over T}}Y(v,x_{1})\Y(w_{(1)},x_{2})w_{(2)}
\nonumber\\
&=&(x_{1}-x_{2})^{k+{j_1\over T}}Y(v,x_{1})y^{L(0)}\psi
\left(y^{-L(0)}w_{(1)}\otimes 
y^{-L(0)}w_{(2)}\right)|_{y=x_{2}e^{-\log z}}
\nonumber\\
&=&(x_{1}-x_{2})^{k+{j_1\over T}}y^{L(0)}Y_W^{*}(y^{-L(0)}v,x_{1}y^{-1})
\psi\left(y^{-L(0)}w_{(1)}\otimes 
y^{-L(0)}w_{(2)}\right)|_{y=x_{2}e^{-\log z}}\nonumber\\
&=&e^{(k+{j_1\over T})\pi i}(x_{2}-x_{1})^{k+{j_1\over T}}y^{L(0)}
Y_{\sigma}^{R}(y^{-L(0)}v,x_{1}y^{-1})
\psi\left(y^{-L(0)}w_{(1)}
\otimes y^{-L(0)}w_{(2)}\right)|_{y=x_{2}e^{-\log z}}
\nonumber\\
&=&e^{(k+{j_1\over T})\pi i}(x_{2}-x_{1})^{k+{j_1\over T}}y^{L(0)}
\psi\left(y^{-L(0)}w_{(1)}\otimes 
Y(y^{-L(0)}v,x_{1}y^{-1})
y^{-L(0)}w_{(2)}\right)|_{y=x_{2}e^{-\log z}}
\nonumber\\
&=&e^{(k+{j_1\over T})\pi i}(x_{2}-x_{1})^{k+{j_1\over T}}y^{L(0)}
\psi\left(y^{-L(0)}w_{(1)}\otimes 
y^{-L(0)}Y(v,x_{1})w_{(2)}\right)|_{y=x_{2}e^{-\log z}}\nonumber\\
&=&e^{(k+{j_1\over T})\pi i}(x_{2}-x_{1})^{k+{j_1\over T}}
\Y(w_{(1)},x_{2})Y(v,x_{1})w_{(2)}.
\end{eqnarray*}

Similarly, for generalized weak associativity relation, let $l\in {\N}$ be such that
$$(x+z)^{l+{j_2\over T}}Y_W^{*}(v,x+z)\psi(w_{(1)}\otimes w_{(2)})
=(z+x)^{l+{j_2\over T}}Y_{\sigma}^{L}(v,x)\psi(w_{(1)}\otimes w_{(2)}).$$
Then 
\begin{eqnarray*}
& &(zx_{0}x_{2}^{-1}+z)^{l+{j_2\over T}}
Y_W^{*}(y^{-L(0)}v,zx_{0}x_{2}^{-1}+z)
\psi(y^{-L(0)}w_{(1)}\otimes y^{-L(0)}w_{(2)})\nonumber\\
&=&(z+zx_{0}x_{2}^{-1})^{l}Y_{\sigma}^{L}(y^{-L(0)}v,zx_{0}x_{2}^{-1})
\psi(y^{-L(0)}w_{(1)}\otimes y^{-L(0)}w_{(2)}).
\end{eqnarray*}
Multiplying both sides by $z^{-l-{j_2\over T}}x_{2}^{l+{j_2\over T}}$ we get
\begin{eqnarray}
& &(x_{0}+x_{2})^{l+{j_2\over T}}Y_W^{*}(y^{-L(0)}v,zx_{0}x_{2}^{-1}+z)
\psi(y^{-L(0)}w_{(1)}\otimes y^{-L(0)}w_{(2)})\nonumber\\
&=&(x_{2}+x_{0})^{l+{j_2\over T}}Y_{\sigma}^{L}(y^{-L(0)}v,zx_{0}x_{2}^{-1})
\psi(y^{-L(0)}w_{(1)}\otimes y^{-L(0)}w_{(2)}).
\end{eqnarray}
Then we obtain the generalized weak associativity relation as
\begin{eqnarray*}
& &(x_{0}+x_{2})^{l+{j_2\over T}}
Y(v,x_{0}+x_{2})\Y(w_{(1)},x_{2})w_{(2)}\nonumber\\
&=&(x_{0}+x_{2})^{l+{j_2\over T}}Y(v,x_{0}+x_{2})y^{L(0)}
\psi\left(y^{-L(0)}w_{(1)}\otimes 
y^{-L(0)}w_{(2)}\right)|_{y=x_{2}e^{-\log z}}\nonumber\\
&=&(x_{0}+x_{2})^{l+{j_2\over T}}y^{L(0)}
Y_W^{*}(y^{-L(0)}v,x_{0}y^{-1}+x_{2}y^{-1})
\psi\left(y^{-L(0)}w_{(1)}\otimes 
y^{-L(0)}w_{(2)}\right)|_{y=x_{2}e^{-\log z}}
\nonumber\\
&=&(x_{0}+x_{2})^{l+{j_2\over T}}y^{L(0)}
Y_W^{*}(y^{-L(0)}v,zx_{0}x_{2}^{-1}+z)
\psi\left(y^{-L(0)}w_{(1)}\otimes 
y^{-L(0)}w_{(2)}\right)|_{y=x_{2}e^{-\log z}}
\nonumber\\
&=&(x_{2}+x_{0})^{l+{j_2\over T}}y^{L(0)}
Y_{\sigma}^{L}(y^{-L(0)}v,zx_{0}x_{2}^{-1})
\psi\left(y^{-L(0)}w_{(1)}\otimes 
y^{-L(0)}w_{(2)}\right)|_{y=x_{2}e^{-\log z}}\nonumber\\
&=&(x_{2}+x_{0})^{l+{j_2\over T}}y^{L(0)}
\psi\left(Y(y^{-L(0)}v,zx_{0}x_{2}^{-1})
y^{-L(0)}w_{(1)}\otimes 
y^{-L(0)}w_{(2)}\right)|_{y=x_{2}e^{-\log z}}\nonumber\\
&=&(x_{2}+x_{0})^{l+{j_2\over T}}y^{L(0)}
\psi\left(y^{-L(0)}Y(v,x_{0})w_{(1)}\otimes 
y^{-L(0)}w_{(2)}\right)|_{y=x_{2}e^{-\log z}}\nonumber\\
&=&(x_{2}+x_{0})^{l+{j_2\over T}}\Y(Y(v,x_{0})w_{(1)},x_{2})w_{(2)}.
\end{eqnarray*}
Then by Lemma \ref{rtwistedcommassoc} and Remark \ref{rintL(-1)}, 
$\Y$ is an intertwining operator. It is clear that $F^{(z)}(\Y)=F_{\Y}^{(z)}=\psi$.
This completes the proof.
\end{proof}

Using essentially the same arguments of \cite{hl2} (cf. \cite{fhl}) we have:

\bp{phl2} 
Let $g_1,g_2,g_3$ be mutually commuting finite-order automorphisms of a vertex operator algebra $V$
with $g_3=g_1g_2$. Let $M_r$ be a $g_r$-twisted $V$-module for $r=1,2,3$. Then
\begin{eqnarray}
I{M_{3}\choose M_{1}\; M_{2}}\cong I{M_{3}\choose M_{2}\; M_{1}},\   \   \   \
I{W_{3}\choose M_{1}\; M_{2}}\cong I{M_{2}'\choose M_{1}\; M_{3}'}
\end{eqnarray}
as vector spaces.
\ep

\br{r-schur}
{\em Since homogeneous subspaces of twisted $V$-modules are finite-dimensional by assumption and we work on ${\C}$, 
for any irreducible twisted module Schur's lemma holds (cf. \cite{fhl}, Remark 4.7.1).
Let ${\rm Irr}_{\sigma_r}(V)$ be a complete set of 
equivalence class representatives of irreducible $\sigma_r$-twisted $V$-modules for $r=1,2$. 
From \cite{fhl}, $W_{1}\otimes W_{2}$ for $W_r\in {\rm Irr}_{\sigma_r}(V)$, $r=1,2$, form a
complete set of equivalence class representatives of irreducible 
$\sigma_1\otimes \sigma_2$-twisted $V\otimes V$-modules.}
\er

Now, we assume $\sigma_2=\sigma_1^{-1}$ and $\sigma=1$.
Let $W_{1}$ and $W_{2}$ be $\sigma_1$-twisted $V$-modules.
The same argument in \cite{li-form} (Remark 2.9) shows that
$I{W_{2}\choose V\; W_{1}}\cong {\rm Hom}_{V}(W_{1},W_{2})$ as vector spaces.
If both $W_{1}$ and $W_{2}$ are irreducible, then by Schur's lemma
we have $\dim I{W_{2}\choose V\; W_{1}}=1$ when $W_{1}\cong W_{2}$, and $0$ otherwise.
Combining this with Proposition \ref{phl2}, we conclude that
$\dim I{V'\choose W_{1}\; W_{2}'}=1$ if $W_{1}\cong W_{2}$, and $\dim I{V'\choose W_{1}\; W_{2}'}=0$ otherwise.
To summarize we have:

\bt{twisted-PW}
Let $\sigma_1$ be a finite-order automorphism of a vertex operator algebra $V$ and let $z$ be a nonzero complex number. 
Then 
\begin{eqnarray}
{\rm soc}\left({\mathfrak{D}}_{\sigma_1,\sigma_1^{-1}}^{(z)}(V)\right)
= \coprod_{W_1\in {\rm Irr}_{\sigma_1}(V)}W_1\otimes W_1'
\end{eqnarray}
as weak $\sigma_1\otimes \sigma_1^{-1}$-twisted $V\otimes V$-modules, where 
${\rm Irr}_{\sigma_1}(V)$ is any complete set of 
equivalence class representatives of irreducible $\sigma_1$-twisted $V$-modules.
\et

For the rest of this section, assume that $\sigma$ is a finite-order automorphism of $V$ with a positive integer $T$ 
such that $\sigma^T=1$.
Let $(W,Y_W)$ be an (ordinary) $\sigma$-twisted $V$-module. 
For $v\in V,\  w\in W$, define
\begin{eqnarray}
\Y_{WV}^{W}(w,x)v=e^{xL(-1)}Y_{W}(v,e^{\pi i}x)w.
\end{eqnarray}
From \cite{fhl}, $\Y_{WV}^W(\cdot,x)$ is an intertwining operator of type $\binom{W}{W\; V}$.
Furthermore, define a linear map $\Y_{W W'}^{V'}(\cdot,x)$ by 
\begin{eqnarray}
\< \Y_{W W'}^{V'}(w,x)w',v\>
&=&\<w', \Y_{W V}^{W}(e^{xL(1)}e^{\pi iL(0)}x^{-2L(0)}w,x^{-1})v\>\nonumber\\
&=&\<w', e^{x^{-1}L(-1)}Y_{W}(v,e^{\pi i}x^{-1})e^{xL(1)}e^{\pi iL(0)}x^{-2L(0)}w\>.
\end{eqnarray}
Then $\Y_{WW'}^{V'}(\cdot,x)$ is an intertwining operator of the indicated type, 
which we call the {\em standard intertwining operator}.
Correspondingly, we have {\em the standard $P(z)$-intertwining map} of type $\binom{V^{*}}{W\;W'}$:
\begin{eqnarray}
\mathcal{F}_W^{(z)}: && W\otimes W'\rightarrow V^{*}\nonumber\\
&&\mathcal{F}_W^{(z)}(w\otimes w')=\Y_{WW'}^{V'}(w,e^{\log z})w'.
\end{eqnarray}

For a $\sigma$-twisted $V$-module $W=\bigoplus_{h\in \C}W_{(h)}$, define a $q$-graded trace function by
\begin{eqnarray}
\tilde{\chi}_W(v,q)=\sum_{h\in \C}{\rm tr}_{W_{(h)}}\Res_x x^{-1}Y_W(x^{L(0)}v,x)q^{L(0)}=\sum_{h\in \C}q^{h}{\rm tr}_{W_{(h)}}v^{W}_{\wt v-1}
\end{eqnarray}
for homogeneous vector $v\in V$, where $Y_W(v,x)=\sum_{n\in \frac{1}{T}\Z}v_n^W x^{-n-1}$.
 Note that if $v\in V^j$ with $0<j<T$, then $v^W_{m}=0$ on $W$ for $m\in \Z$,
 in particular, $\wt v^W_{\wt v-1}=0$ (as $\wt v\in \Z$). Thus
 \begin{eqnarray}
 \tilde{\chi}_W(v,q)=0\quad \text{for }v\in V^j,\ 0<j<T.
 \end{eqnarray}

If $W$ is indecomposable, there exists a unique $\lambda\in \C$ such that
$W_{(\lambda)}\ne 0$ and $W=\oplus_{n\in \N}W_{(\lambda+\frac{n}{T})}$. 
Then $\tilde{\chi}_W(v,q)\in q^{\lambda}\C[[q^{\frac{1}{T}}]]$ for $v\in V$. That is,
\begin{eqnarray}
\tilde{\chi}_W(\cdot,q)\in q^{\lambda}V^{*}[[q^{\frac{1}{T}}]].
\end{eqnarray}

\bp{thm-trace-function}
Let $W$ be a $\sigma$-twisted $V$-module. Then the coefficients of $q^{h}$ for $h\in \C$ in
the $q$-graded trace function $\tilde{\chi}_W(\cdot,q)$ lie in $\mathfrak{D}_{\sigma,\sigma^{-1}}^{(-1)}(V)$.
Furthermore, if $W$ is irreducible, then all the coefficients in $\tilde{\chi}_W(\cdot,q)$ 
generate an irreducible $\sigma\otimes \sigma^{-1}$-twisted $V\otimes V$-submodule
of $\mathfrak{D}_{\sigma,\sigma^{-1}}^{(-1)}(V)$, which is isomorphic to $W\otimes W'$. 
\ep

\begin{proof} 
Let $\mathcal{B}=\{ w_{\alpha}\ |\ \alpha\in S\}$ be a basis of $W$ with all $w_{\alpha}$ homogeneous,
and let $\mathcal{B}'=\{ w^{\alpha}\ |\ \alpha\in S\}$ be the corresponding dual basis of $W'$.
Let $v\in V^{0}$ homogeneous. Then
\begin{eqnarray}
&&\tilde{\chi}_W(v,q)\nonumber\\
&&=\sum_{\alpha\in S}\<w^{\alpha},v_{\wt v-1}q^{L(0)}w_{\alpha}\>\nonumber\\
&&=\sum_{\alpha\in S}\<w^{\alpha},Y_W((e^{\pi i}x)^{L(0)}v,e^{\pi i}x)q^{L(0)}w_{\alpha}\>\nonumber\\
&&=\sum_{\alpha\in S}\<w^{\alpha},e^{-xL(-1)}\Y_{WV}^W(q^{L(0)}w_{\alpha},x)(e^{\pi i}x)^{L(0)}v\>\nonumber\\
&&=\sum_{\alpha\in S}\<e^{-xL(1)}w^{\alpha},\Y_{WV}^W(q^{L(0)}w_{\alpha},x)(e^{\pi i}x)^{L(0)}v\>\nonumber\\
&&=\sum_{\alpha\in S}\<\Y_{WW'}^{V'}\left(x^{-2L(0)}e^{-\pi i L(0)}e^{-x^{-1}L(1)}q^{L(0)}w_{\alpha},x^{-1}\right)
e^{-xL(1)}w^{\alpha},(e^{\pi i}x)^{L(0)}v\>\nonumber\\
&&=\sum_{\alpha\in S}\<e^{\pi iL(0)}x^{L(0)}\Y_{WW'}^{V'}\left(x^{-2L(0)}e^{-\pi i L(0)}e^{-x^{-1}L(1)}q^{L(0)}w_{\alpha},x^{-1}\right)e^{-xL(1)}w^{\alpha},v\>\nonumber\\
&&=\sum_{\alpha\in S}\<\Y_{WW'}^{V'}\left(e^{\pi iL(0)}x^{-L(0)}e^{-\pi i L(0)}e^{-x^{-1}L(1)}q^{L(0)}w_{\alpha},-1\right)
e^{\pi iL(0)}x^{L(0)}e^{-xL(1)}w^{\alpha},v\>\nonumber\\
&&=\sum_{\alpha\in S}\<\Y_{WW'}^{V'}\left(x^{-L(0)}e^{-x^{-1}L(1)}q^{L(0)}w_{\alpha},-1\right)e^{\pi iL(0)}x^{L(0)}e^{-xL(1)}w^{\alpha},v\>\nonumber\\
&&=\sum_{\alpha\in S}\<\Y_{WW'}^{V'}\left(e^{-L(1)}q^{L(0)}x^{-L(0)}w_{\alpha},-1\right)e^{\pi iL(0)}e^{-L(1)}x^{L(0)}w^{\alpha},v\>\nonumber\\
&&=\sum_{\alpha\in S}\<\Y_{WW'}^{V'}\left(e^{-L(1)}q^{L(0)}w_{\alpha},-1\right)e^{\pi iL(0)}e^{-L(1)}w^{\alpha},v\>\nonumber\\
&&=\sum_{\alpha\in S}\<\Y_{WW'}^{V'}\left(e^{-L(1)}q^{L(0)}w_{\alpha},-1\right)e^{L(1)}e^{\pi iL(0)}w^{\alpha},v\>\nonumber\\
&&=\sum_{\alpha\in S}\<\Y_{WW'}^{V'}\left(e^{-L(1)}e^{\pi iL(0)}q^{L(0)}w_{\alpha},-1\right)e^{L(1)}w^{\alpha},v\>,
\end{eqnarray}
noticing that $\wt w_{\alpha}=\wt w^{\alpha}$ for $\alpha\in S$ and
$$x_1^{-L(0)}e^{xL(1)}x_1^{L(0)}=e^{xx_1L(1)}.$$
In view of Theorem \ref{t-Int-map-dual},  the standard $P(-1)$-intertwining map $\mathcal{F}_W^{(-1)}$ 
is a homomorphism of weak $\sigma\otimes \sigma^{-1}$-twisted 
$V\otimes V$-modules from $W\otimes W'$ into $\mathfrak{D}_{\sigma,\sigma^{-1}}^{(-1)}(V)$.
Then the first assertion follows immediately. The second assertion is clear as $W\otimes W'$ is an 
irreducible $\sigma\otimes \sigma^{-1}$-twisted 
$V\otimes V$-module and $\mathcal{F}_W^{(-1)}$ is a nontrivial module homomorphism.
\end{proof}

\end{document}